\documentclass{article}

\usepackage{amsfonts}
\usepackage{color}
\usepackage{graphics}
\usepackage{epsfig,psfrag,graphicx}
\usepackage{amssymb}
\usepackage{eepic,epic}
\usepackage{epsfig} 
\usepackage{amsmath} 
\usepackage{amssymb} 
\usepackage{amsbsy} 
\usepackage{fullpage}

\def\XXint#1#2#3{{\setbox0=\hbox{$#1{#2#3}{%
\int}$ }
\vcenter{\hbox{$#2#3$ }}\kern-.6\wd0}}

\usepackage{amsmath}

  \usepackage{paralist}
  \usepackage{graphics} 
  \usepackage{epsfig} 
 \usepackage[colorlinks=true]{hyperref}
\hypersetup{urlcolor=blue, citecolor=red}

\usepackage{amssymb}
\usepackage{amsfonts}
\usepackage{color}
\usepackage{graphics}
\usepackage{epsfig,psfrag,graphicx}
\usepackage{amssymb}
\usepackage{eepic,epic}
\usepackage{epsfig} 
\usepackage{amsmath} 
\usepackage{amssymb} 
\usepackage{amsbsy} 

\newtheorem{theorem}{Theorem}[section]

\newtheorem{lemma}[theorem]{Lemma}

\newtheorem{definition}[theorem]{Definition}

\newcommand{\Div}       {{\rm div}_x}

\newcommand{\D}     {{\mathcal D}}

\def\tens#1{\pmb{\mathsf{#1}}}
\def\vec#1{\boldsymbol{#1}}

\newcommand{\uu}        {{\vec u}}
\newcommand{\uke}        {\uu^{k,\epsilon} }
 \newcommand{\uk}        {\uu^{k} }
 \newcommand{\un}        {\uu^{n} }

\newcommand{\Grad}      {\nabla_x}
\newcommand{\vr}        {\varrho}
\newcommand{\vrke}        {\varrho^{k,\epsilon}}
\newcommand{\vrk}        {\varrho^{k}}
\newcommand{\vrn}        {\varrho^{n}}

\newcommand{\tke}        {\theta^{k,\epsilon} }
\newcommand{\tk}        {\theta^k }
\newcommand{\tn}        {\theta^n }
\newcommand{\ake}        {\alpha^{k,\epsilon} }
\newcommand{\bke}        {\nu^{k,\epsilon} }

\newcommand{\V}         {\mathcal{V}}
\newcommand{\R}         {\mathbb{R}}
\newcommand{\N}         {\mathbb{N}}

\newcommand{\grs}				{\tens{D}}

\newcommand{\E}		{\vec{E}}
\newcommand{\Du}        {\tens{D}\uu}
\newcommand{\Dun}				{\tens{D}\uu^n}
\newcommand{\bS}        {\tens{S}}
\newcommand{\Sdu}       {\tens{S}(x,\vr, \theta, \tens{D}\uu)}

\newcommand{\Sdurt}      {\tens{S}(x,\vr^n,\tn, \tens{D}\uu)}

\newcommand{\Ssa}       {\tens{S}(x,l, s, \lambda)}

\newcommand{\Ssl}       {\tens{S}(x,\vr, \theta, \lambda)}

\newcommand{\Sk}       {\tens{S}_{k} }

\newcommand{\Ske}       {\tens{S}_{k,\epsilon} }
\newcommand{\qu}       {\vec{q} (\vr, \theta , \Grad \theta)}
\newcommand{\qn}       {\vec{q} (\vrn, \tn , \Grad \tn)}

\newcommand{\Sdun}      {\tens{S}(x,\vr^n,\tn, \tens{D}\uu^n)}

\newcommand{\bff}       {{\vec{f}}}

\newcommand{\oo}				{{\vec{\omega}}}

\newcommand{\bchi}			{\overline{\bS}}

\newcommand{\Vpd}				{W^{1,p}_{0,{\rm div}}(\Omega)^3}
\newcommand{\Vdpd}				{W^{1,2p}_{0,{\rm div}}(\Omega)^3}

\newcommand{\cp}        {\underline{c}}

\newcommand{\vd}        {{\tens{w}}}
\newcommand{\bvarphi}   {{\vec{\varphi}}}
\newcommand{\zd}        {{\tens{v}}}

\newcommand{\bbK}         {\tens{K}}
\newcommand{\bbL}         {\tens{L}}

\newcommand{\txi}		{\tens{\xi}}

\newcommand{\tz}		{\tens{z}}
\newcommand{\weak}			{\rightharpoonup}
\newcommand{\weakstar}	{\stackrel{\ast}{\rightharpoonup}}
\newcommand{\modular}[1]{{\stackrel{ #1}{\longrightarrow\,}}}
\newcommand{\weakly}			{\quad \mbox{ weakly in }}
\newcommand{\biting}          {\rightarrow_b}
\def\bbbone{{\mathchoice {\rm 1\mskip-4mu l}
{\rm 1\mskip-4mu l} {\rm 1\mskip-4.5mu l} {\rm 1\mskip-5mu l}}}

\title
      {Unsteady flows of  heat-conducting non-Newtonian  fluids in~Musielak-Orlicz spaces}

\author{
Bart\l omiej Matejczyk\footnote{
Johann Radon Institute for Computational and Applied Mathematics (RICAM),
Austrian Academy of Sciences, Altenberger Stra{\ss}Ÿe 69, 4040 Linz, Austria, bartlomiej.matejczyk@oeaw.ac.at
} 
\ and  Aneta Wr\'oblewska-Kami\'nska\footnote{
Institute of Mathematics of the Polish Academy of Sciences, ul. \'Sniadeckich 8, 00-656, Poland, awrob@impan.pl }
}

\begin{document}



%

\maketitle



%
%
%
%

\bigskip


\begin{abstract}
Our purpose is to show the existence of weak solutions for unsteady flow of non-Newtonian incompressible nonhomogeneous, heat-conducting fluids with  generalised form of the  stress tensor without any restriction on its upper growth. Motivated by fluids of nonstandard rheology we focus on  the general form of growth conditions for the stress tensor which makes anisotropic Musielak-Orlicz spaces a suitable function space for the considered problem. We do not assume any smallness condition on initial data in order to obtain long-time existence. Within the proof we use monotonicity methods, integration by parts adapted to nonreflexive spaces and Young measure techniques.\\
{{\bf Subclass:} 35Q35, 46E30, 76D03, 35D30}

\end{abstract}

\section{Introduction and formulation of the problem}\label{introduction}
In this article we investigate mathematical model of the flow of an  incompressible, nonhomogeneous non-Newtonian, heat-conducting  fluid governed by the following system of equations:
	\begin{equation}\label{PP}
        \begin{split}
        \partial_t\vr + \Div (\vr \uu)&=0 \quad{\rm in}\quad Q,\\
        \partial_t(\vr\uu) + \Div(\vr\uu\otimes\uu) - \Div\Sdu + \Grad P &= \vr \bff  \quad{\rm in}\quad Q,\\
     \partial_t(\vr\theta) + \Div(\vr\uu\theta) - \Div\qu  &=\Sdu  : \Du \quad{\rm in}\quad Q,\\
        \Div\uu&=0 \quad{\rm in}\quad Q,\\
        \uu(0,x)&=\uu_0\,\,\,\,{\rm in}\quad\Omega,\\
        \vr(0,x)&=\vr_0\,\,\,\,{\rm in}\quad\Omega,\\
       \vec q  \cdot \vec{n}   = 0, \quad  \uu(t,x)&=0\quad{\rm on}\quad[0,T]\times\partial\Omega,
        \end{split}
        \end{equation}
where $\vr:Q\to\R$ is a mass density, $\uu:Q\to\R^3$ stands for  a velocity field, $P:Q\to\R$ is  a 
pressure function, $\bS$ -  a stress tensor, $\vec{q}$ - a thermal flux vector, $\bff :Q\to\R^3$ - a given outer force. The set $\Omega\subset\R^3$
is a bounded domain with a regular boundary $\partial\Omega$ (of class, 
say $C^{2+\nu}$, $\nu>0$, taken for convenience). We
denote by $Q=(0,T)\times\Omega$ the time-space cylinder with some given $T\in(0,+\infty)$. The tensor 
$\Du=\frac{1}{2}(\Grad \uu+\Grad^T\uu)$ is a symmetric part of the velocity gradient.

For the above system we set the initial density $\vr$  and temperature $\theta$  to satisfy  
  \begin{equation}  \label{as-vr0}
    \vr(0,\cdot)=\vr_0\in L^{\infty}(\Omega)  \quad \mbox{ and }  \quad 
    0<\vr_*\leq\vr_0(x)\leq\vr^*<+\infty\quad{\rm for\,a.a.\,}x\in\Omega , 
    \end{equation}
    
 \begin{equation}
 \label{as-th0}\theta_0 \in L^1(\Omega) \quad \mbox { and } \quad 
    0<\theta_*\leq\theta_0 \quad{\rm for\,a.a.\,}x\in\Omega.
    \end{equation}

In order to formulate  the growth conditions of the stress tensor we use 
general convex function $M$ called an $N$--function  similarly as in  \cite{gwiazda-swierczewska,GSW2,AW,AW2013} (for a definition see Section \ref{orlicz}). We assume that stress tensor 
 $\bS: \Omega\times \R_+ \times \R_+ \times \R^{3\times 3}_{\rm{sym}}\to\R^{3\times
3}_{\rm{sym}}$ satisfies ($\R^{3\times 3}_{\rm{sym}}$ stands for
the space of $3\times 3$ symmetric matrices):

\begin{description}
    \item[\bf S1.] $\bS(x,\vr,\theta,{\bbK}) $ is a Carath\'eodory function (i.e., measurable
    function of $x$ for all $\vr, \theta > 0$ and ${\bbK}\in\R^{3\times 3}_{\rm{sym}}$ and
    continuous function of $\theta, \vr$ and ${\bbK}$ for a.a. $x\in\Omega$) and
    $\bS(x,\vr,\theta,\tens{0})=\tens{0}$.

    \item[\bf S2.] There exists a positive constant $c_c \in (0,1)$, $N$--functions $M$
    and $M^\ast$ (which denotes the complementary function to $M$)  such that for all $\bbK\in\R^{3\times 3}_{\rm{sym}}$, $\theta, \vr>0$  and
    a.a. $t,\ x\in Q$  holds
        \begin{equation}\label{coercivity}
        \bS(x,\vr, \theta,{\bbK}) : {\bbK}\geq
        c_c\{M(x,{\bbK})+M^*(x,\bS(x,\vr,\theta,{\bbK}) )\}.
        \end{equation}


    \item[\bf S3.] $\bS$ is monotone, i.e. for all ${\bbK}_1,{\bbK}_2\in\R^{3\times
    3}_{\rm{sym}},$ $\vr>0$, $\theta >0$ and a.a. $x\in\Omega$
    $$[\bS(x,\vr,\theta,{\bbK}_1)-\bS(x,\vr,\theta,{\bbK}_2)] : [{\bbK}_1-{\bbK}_2]\geq0.$$
\end{description}

The heat flux  $\vec q $, as usually,  is set to be less general.  Therefore, similarly as in  \cite{FMR},  we expect  $\vec q $ of the form
\begin{equation*}
\qu \approx \kappa (\vr) \theta^\beta \Grad \theta = \kappa(\vr) \frac{1}{\beta+1} \nabla \theta^{\beta + 1}  \quad \mbox{ for } \beta \in \R,  
\end{equation*}
such that $\kappa (\vr)$ satisfies 
$
0\leq \kappa_* \leq \kappa(\vr) \leq \kappa^* < \infty , 
$
where $\kappa_*, \kappa^*$ are some fixed constants.  In particular, we require that $\vec{q}: \R_+ \times \R_+ \times \R^3 \to \R^3$ satisfies
	\begin{equation}
	\vec{q}(\vr,\theta,\nabla_x\theta) = \kappa_0 (\vr,\theta) \nabla_x \theta  \quad  \mbox { with   } \kappa_0 \in C(\R^2_+)
	\label{zal_q}
	\end{equation} 
and 
 for all $ \theta , \vr  > 0 ,\ \Grad \theta \in \R^3  $
\begin{equation}\label{ass-q}
\begin{split}
\qu  \cdot \Grad \theta & \geq
 \kappa_* \theta^\beta |\Grad \theta|^2  \quad  \mbox{ with } \beta \in \R \mbox{ and } \kappa_* >0  , \\
 |\qu| & \leq \kappa^* \theta^\beta | \nabla \theta|  \quad  \quad \mbox{ with } \quad   \kappa^* >0 .
 \end{split}
\end{equation}

\bigskip

		 The  main  reason to investigate such a  general form of stress tensor $\bS$, namely satisfying  \eqref{coercivity}, is the phenomena of rapidly increasing fluid viscosity under various stimuli as shear rate,  electric or magnetic field. 
		Our assumptions include power-law and 
Carreau-type models which are quite popular among rheologists, 
chemical engineering and colloidal mechanics.

		 In majority of publications concerning non-Newtonian fluids a $p$-structure for $\bS$ is assumed. It means that 
$\bS\approx \mu(\vr,\theta)(\kappa + |\Du|)^{p-2}\Du$ or 
$\bS\approx \mu(\vr,\theta)(\kappa + |\Du|^2)^{(p-2)/2}\Du$ (where $\kappa >0$ and $\mu$ is a nonnegative bounded function).
Then standard growth conditions of the stress tensor, namely
polynomial growth:
  $
    |\bS(x,\tens\xi)|\leq c(1+|\tens\xi|^2)^{(p-2)/2}|\tens\xi|$ and 
    $
    \bS(x,\tens\xi):\tens\xi \geq c(1+|\tens\xi|^2)^{(p-2)/2}|\tens\xi|^2
$ are satisfied, see e.g. \cite{FMR,frehse,MaNeRu93}. Unfortunately this theory is  not adequate for phenomena of fluids that rapidly and significantly  change their viscosity, i.e. when growth of the stress tensor may be much faster then polynomial and which may differ in various directions of shear stress  or be inhomogeneous is spatial variables. Examples of these fluids are shear thickening (STF), magnetorheological (MR) and electorheological fluids.  Because of  property of the changeable  viscosity these fluids have applications  in a variety of industry, military and natural sciences. 

Firstly we would like to be able to consider flows for which constitutive relations for the stress tensor $\bS$ are more general than of power-law type and, in particular, which the growth w.r.t. the shear stress may be faster than polynomial. Very  promising application of this type of fluids  is the one coming from military industry. STF fluids behave as a liquid  until  another object strike it with high kinetic energy. In this case the fluid increases its viscosity in milliseconds and behaves almost like a solid. Moreover this process is completely reversible which  makes such a fluid a perfect material for  military, medical and sport armours.   Obtained material has  high elasticity   combined with protection   against needles, knifes and bullets  \cite{decker2007stab,ELK,HSK,YLWW}.

Moreover,  we can study constitutive relations for fluids with dependence on outer field (magnetic or electric). 
 Mathematical models for such fluids are  considered e.g. in \cite{RaRu2000}. 
Governing equations are derived from motion of electrorheological fluids,  taking into consideration complex interactions between electromagnetic  and thermomechanical fields  into consideration (see also \cite{RaRu1996}). 
For such general fluids, as claimed in (cf. \cite{Ruzicka}),  the  stress tensor can be written in a quite general form, which is  still thermodynamically admissible (i.e. $\bS : \grs \geq 0$), satisfies the principle of material frame-indifference and is monotone. But then it may appear  that the standard growth conditions, i.e. $|\bS(\grs,\E) | \leq c(1 + |\grs | )^{p-1},$ $\bS(\grs,\E): \grs \geq c |\grs|^p$
($\E$ denotes electric flux) are not satisfied, because the tensor $\bS$ may possess the  growth of different powers in various directions of $\grs$ (for the example see also \cite{AW2013}).

In our considerations we also would like to cover the case of constitutive relations which may depend on spatial variables. For example, again it may be  the  case of  electrorheological  fluids  which are suspensions of extremely fine non-conducting particles  in an electrically insulating fluid. Such a mode was considered e.g. in \cite{Ruzicka} where  the $N$--function took the form: $M (x,z) = |z|^{p(x)}$ with $1< p_{-} \leq p(x)\leq p+ <\infty$.  The author provided there the existence theory  for the case of barotropic flows without dependence on density.

\medskip

 The appropriate spaces to capture such formulated problem are anisotropic Musielak-Orlicz spaces.
 For definitions and preliminaries of $N$--functions and Orlicz spaces see Section \ref{orlicz}.
 Contrary to \cite{Musielak} we consider the $N$--function $M$ not
 dependent only on $|\tens\xi|$, but on whole tensor $\tens\xi.$
 It results from the fact that the viscosity may differ in different directions of symmetric part of velocity gradient
 $\grs \uu$ and 
 the growth condition for the stress tensor dependent  on the whole tensor $\grs \uu$, not only on $|\grs \uu|$.
 Since we allow $\bS$ to depend on spatial variable, $N-$function depends also on $x \in \Omega.$
 The general growth on $\bS$ is provided by quite general properties of the $N-$function $M$ defining anisotropic Musielak-Orlicz space $L_M$. Since we do not want to be restricted by any upper growth conditions on $\bS$, we do not assume that, so called, $\Delta_2-$condition is satisfied by $M$. 
 The spaces with $N$--function dependent on vector-valued argument were introduced in \cite{skaff1,skaff2,Turett}.
Let us underline that, in general, if $M$ and $M^*$  do not satisfy a $\Delta_2$-condition the related spaces fail to be separable and reflexive, which is a source  of additional difficulties arising from functional analysis. We then lose simply lose many facilitating 
properties of spaces we have work with. 
The setting considered in this paper needs 
tools which generalise results not only of classical Lebesgue and Sobolev spaces (related to power-law fluids), but also these in variable exponent, anisotropic and classical Orlicz spaces.
Most of the essential and necessary tools of functional analysis for classical Orlicz spaces (isotropic and homogenous) are already deeply understood, for example the density of soothe functions in modular topology \cite{Gossez} and integration by parts formula \cite{EM}. But many structures for anisotropic Musielak-Orlicz spaces have not been developed or are not understood purely yet.

\medskip

One of the essential difficulty we have to face to provide is the weak sequential stability in energy equation. Namely we need to show that	$\bS^n := \bS (\cdot, \varrho^n, \theta^n, \tens{D} \vec{u}^n) : \tens{D} \vec{u}^n \weak \bS (\cdot, \varrho, \theta, \tens{D} \vec{u}) : \tens{D} \vec{u}$ weakly in $L^1(Q)$,  where  $\{\varrho^n\}_{n=1}^\infty ,$ $\{\vec{u}^n\}_{n=1}^\infty ,$ $\{\theta^n \}_{n=1}^\infty $ are approximation sequences  of 
$\varrho,$ $\vec{u},$ $\theta$ and $\{ \bS^n \}_{n=1}^\infty \subset L_{M^*}(Q)$, $\{ \tens{D} \vec{u}^n \}_{n=1}^\infty \subset L_{M}(Q)$. 
Let us notice that if one work with reflexive spaces (such a $L^p$) the monotonicity is a sufficient argument to conclude from $(\bS^n - \bS):(\tens{D}\vec{u}^n - \tens{D}\vec{u})  \weak 0$ in $L^1$ that $\bS^n :  \Dun  \weak \bS :  \Du $ weakly in $L^1$. However, once the space is not reflexive, as the case of our $L_M$-space, then the convergence may fail if one is not able to provide modular convergence of sequences $\bS^n$ and $\Dun$ in proper spaces. 
In the current paper we use bitting lemma \cite{Aliberti,ball,pedregal} and methods of Young measures to show that the product of our two sequences converges weakly in $L^1$ and consequently to provide the sequential stability of the right hand side of energy equation.
Similar arguments in frame of anisotropic Musielak-Orlicz spaces  were used in \cite{gwiazda2015renormalized} for parabolic equation and later also in \cite{Klawe} for the problem of thermo-visco-elasticity model.

An interesting obstacle here is the lack of the classical
integration by parts formula, see~\cite[Section 4.1]{GGZ}. To
extend it for the case of anisotropic Musielak-Orlicz spaces we  would need that $C^\infty-$functions are  dense in $L_M(Q)$
and $L_M(Q)=L_M(0,T;L_M(\Omega))$. The first one only holds if $M$
satisfies $\Delta_2$--condition. The second one is not the case in
Orlicz and generalized Orlicz spaces, see \cite{Donaldson} and holds only if $M$
is equivalent to some power $p$, $1< p<\infty$ (what provides that $L_M(Q)$ is separable and reflexive).
In the present paper we recall the integration by parts formula obtained in \cite{AW2013} by adaptation of arguments from \cite{GSW2} and \cite{FMR}.

Moreover classical monotonicity methods allowing to obtain convergence in a nonlinear viscous term in the momentum equation do not work in case of non-reflexive anisotropic Musielak-Orlicz spaces. Therefore we need to apply arguments developed in \cite{AW,AW2013,GSW2}, see also \cite{MustonenTienari}.

\medskip

Let us now recall briefly related results. 
The mathematical analysis of time dependent flow of homogeneous (density was assumed to be constant) non-Newtonian fluids of power-law type was initiated in \cite{La1,La2}, where the global existence of weak solutions for the exponent $p\geq 1 + (2d)/(d+2)$ ($d$ stands for space dimension) was proved for Dirichlet boundary conditions. 
Later the steady flow was considered in \cite{fms}, where the existence of weak solutions was established for 
the constant exponent $p>\frac{2d}{d+2}$, $d\geq 2$ by Lipschitz truncation methods.

In \cite{Ruzicka} generalized Lebesgue spaces $L^{p(x)}$ were used
 to the description of flow of electrorheological fluid.
 The  author assumed in this work that $1< p_0\leq p(x) \leq p_\infty
<\infty$, where $p\in C^1(\Omega)$ is a function of an electric field $E$, i.e.
 $p=p(|E|^2)$, and $p_0 >\frac{3d}{d+2}$ in case of steady flow, where $d\geq 2$ is the space dimension. The $\Delta_2$--condition is then satisfied  and consequently the space is reflexive
and separable (what is not the case in of our work).  Later in \cite{DMS}  the above result was 
improved by Lipschitz truncation methods for $L^{p(x)}$ setting for $\bS$,
where $\frac{2d}{d+2}<p(\cdot)<\infty$ was log-H\"{o}lder continuous and $\bS$ 
was strictly monotone.

In \cite{WO} the author proved existence  of weak solutions to
unsteady motion of an incompressible homogenous fluid with shear rate dependent viscosity
for $p> 2(d+1)/(d+2)$ without assumptions on shape and size of $\Omega$ employing an $L^\infty$--test function and local pressure method. Finely the existence of global weak solutions  with Dirichlet boundary conditions for $p>(2d)/(d+2)$ was achieved in \cite{drw}  by
Lipschitz truncation and local pressure methods.

Most of the available results concerning
nonhomogeneous (without assumption that density is constant) incompressible fluids deal with the polynomial dependence
between $\bS$ and  $|\Du|$. The analysis of nonhomogeneous Newtonian
($p=2$ in standard growth condition) fluids
was investigated in \cite{kazhikhov} in the seventies.
In \cite{lions}  the concept of renormalized solutions was presented what allowed to obtain convergence and continuity properties of the density.

The first result for unsteady flow of nonhomogenous  non-Newtonian fluids goes back to \cite{FGO}, where existence of Dirichlet weak solutions was obtained  for  $p\geq 12/5$ if $d=3$, later existence of space-periodic weak solutions for $p\geq2$  with some regularity properties of weak solutions whenever $p\geq 20/9$ (if $d=3$)   was achieved  in \cite{GG}.
In \cite{frehse} existence of a weak solution was showed 
for generalized Newtonian fluid of power-law type for $p>11/5$.
Authors also needed existence of the potential of $\bS$.
The most related result concerning inhomogeneous, incompressible and heat-conducting non-Newtonian fluids, but of standard growth conditions of polynomial type 
 for $p\geq 11/5$  the reader can find in  \cite{FMR}. The novelty
of this paper w.r.t. the previously mentioned results was the consideration of the full thermodynamic
model.

The analysis of  non-Newtonian fluids  in frame of anisotropic  Musielak-Orlicz spaces have been studied 
with variety of approaches.  Some of  considerations can be found in \cite{gwiazda-swierczewska} (the case of homogeneous, incompressible non-Newtonian fluids) where $\tens{S}$ was taken strictly monotone. The authors used Young measure technics in place of monotonicity methods. The additional assumption at strict monotonicity allows to conclude that the measure-valued solution is of the form of Dirac measure and then the system has a weak solution. 
Later generalisation of the Browder-Minty trick for non-reflexive anisotropic Musielak-Orlicz was used in   \cite{GSW2,AW,AW2013}, what allows to assume only the monotonicity of  $\tens{S}$.

In \cite{GSW3} authors studied generalized Stokes system for the unsteady flow of homogenous, incompressible non-Newtonian fluids of non-standard rheology. Neglecting the convective term in momentum equation they showed existence of weak solutions in anisotropic Orlicz spaces without assumption on lower bound on $N$--function $M$, what allowed them to consider also shear-thinning fluids.
 
In particular in \cite{AW2013} the author obtained existence of weak solutions to unsteady flow of non-Newtonian incompressible nonhomogeneous fluids with nonstandard growth conditions of the stress tensor assumed also in the current paper.

\medskip

Summarising, our less restrictive assumptions on tensor $\bS$ allow to consider effects of nonhomogeneous (dependence on $x$), anisotropic behaviour of considered medium and as well more general than power-law type rheologies.  
In this  article we focus on time dependent flow of  non-Newtonian, inhomogeneous (density dependent), incompressible fluid and our main goal is to consider also  temperature and its influence on the flow. Let us emphasise that the stress tensor $\bS$ may depend here not only on a shear stress but also on density and changes in temperature. Up to our knowledge, the existing result for such a problem has not been considered yet and our considerations extend the   theory concerning  thermodynamics of non-Newtonian fluids of power-law type to the case of non-standard growth conditions in Musielak-Orlicz spaces on the one hand, and the theory of non-Newtonian fluids with non-standard growth conditions to the case of heat-conducting fluids on the other.

\medskip

	Our paper is constructed as follows in Section \ref{pre1} we recall some used facts and definition necessary for the main theorem stated in Section \ref{defsol}. In  Section~~\ref{proof} we prove the main theorem building $n$--approximate solutions, providing uniform estimates and using monotonicity method, compensated compactness arguments and Young measures.

\section{Preliminaries}\label{pre1}

\subsection{Used notation}\label{orlicz}
In  the following section we introduce   notation, 
definitions and some important properties   of  Orlicz spaces used in further considerations.   More studies of Orlicz spaces can be found in \cite{Kufner,Musielak,skaff1,skaff2} .

By $\D(\Omega)$   we mean the set of $C^\infty$-functions with compact
support contained in $\Omega$.
 Let ${\V}$ be the set of all functions which belong to ${\mathcal
D}(\Omega)$ and are divergence-free. Moreover, by $L^p, W^{1,p}$ we
denote the standard Lebesgue and Sobolev spaces respectively, by $H$--the
closure of $\V$ w.r.t.  the $\|\cdot\|_{L^2}$ norm and by
$W^{1,p}_{0,{\rm div}}$--the closure of $\V$ w.r.t.  the
$\|\nabla(\cdot)\|_{L^p}$ norm. Let $W^{-1,p'}=(W^{1,p}_0)^*$ and 
$W^{-1,p'}_{\rm div}=(W^{1,p}_{0, {\rm div}})^*$.  By $p'$ we denote a conjugate exponent to $p$, namely $\frac{1}{p}+\frac{1}{p'}=1$.

If $X$  is a Banach space of scalar functions, then $X^3$ or $X^{3\times 3}$
denotes the space of vector-  or tensor-valued functions where each
component belongs to $X$. The symbols $L^p(0,T;X)$ and $C([0,T];X)$
mean the standard Bochner spaces.
Finally, we recall that the Nikolskii space $N^{\alpha,p}(0,T ; X)$
corresponding to the Banach space X and the exponents $\alpha\in (0,1)$ and $p\in[1,\infty]$ is given by
$$N^{\alpha,p}(0,T;X):= \{ f \in L^p(0,T ; X)\ : \ \sup\limits_{0<h<T} h^{-\alpha} \| \tau_h f - f \|_{L^p(0,T-h;X)} < \infty  \} ,$$
where $\tau_h f(t) = f( t + h)$ for a.a. $t \in [0,T-h]$. 
By $(a,b)$ we mean $\int_\Omega a(x)\cdot b(x) dx$ an  inner product of two vector functions or in case    $\int_\Omega a(x) :  b(x) dx$  of  two tensor functions and $\left\langle a, b \right\rangle$ denotes the duality pairing.
\subsection{Orlicz spaces}
\begin{definition}
    Let $\Omega$ be a bounded domain in $\R^3$, a function
    $M:\Omega\times\R^{3\times 3}_{\rm{sym}}\to\R_+$ is said to be
    an \emph{$N$--function} if it satisfies the following conditions
    \begin{enumerate}
        \item $M$  is a Carath\'{e}odory function (measurable w.r.t to the first argument and continuous w.r.t. the  second one)  such that $M(x,\bbK)=0$ if and only if
            $\bbK=0$,
            $M(x,\bbK)=M(x,-\bbK)$ a.e. in $\Omega,$
        \item $M(x,\bbK)$ is a convex function w.r.t. $\bbK,$
        \item
            \begin{equation}\label{limitinfty}
            \lim\limits_{|\bbK|\to0}   \frac{M(x,\bbK)}{|\bbK|}=0
            \quad \mbox{ and } \quad
            \lim\limits_{|\bbK|\to\infty}  \frac{M(x,\bbK)}{|\bbK|}=\infty \quad \mbox { for a.a.  } x\in \Omega .
            \end{equation}
    \end{enumerate}
\end{definition}

\begin{definition} The \emph{complementary function  $M^\ast$} to a function $M$ is
    defined  for $\bbL\in\R^{3\times 3}_{\rm{sym}},\; x\in\Omega$ by
    $$M^\ast(x,\bbL)=\sup\limits_{\bbK\in\R^{3\times 3}_{\rm{sym}}}\left(\bbK : \bbL-M(x,\bbK)\right).$$
   
\end{definition}
 Let us notice that the complementary function $M^\ast$ is also an $N$-function (see  \cite{skaff1}). 
\begin{definition} Let $Q = (0,T) \times \Omega$. 
 The \emph{anisotropic Musielak-Orlicz class} ${\mathcal L}_M(Q)^{3\times 3}_{\rm{sym}}$
is the set of all measurable functions
 $\bbK:Q\to\R^{3\times 3}_{\rm{sym}}$ such that
    $$\int_Q M(x,\bbK(t,x)) dxdt<\infty.$$
\end{definition}
\begin{definition}
The \emph{anisotropic Musielak-Orlicz space} $L_M(Q)^{3\times
3}_{\rm{sym}}$ is defined as the set of all measurable functions
$\bbK:Q\to\R^{3\times 3}_{\rm{sym}}$ which satisfy
    $$\int_Q M(x,\lambda \bbK(t,x))dxdt\to0\quad {\rm as}\,\lambda\to
     0.$$
\end{definition}
A Musielak-Orlicz  space is a Banach space with  Luxemburg norm given by
    $$\|\bbK\|_M=\inf\left\{\lambda>0:\quad\int_Q
     M\left(x,\frac{\bbK(t,x)}{\lambda}\right)dxdt\le 1\right\}.$$
Let us denote by $E_M(Q)^{3\times 3}_{\rm{sym}}$ the closure of all measurable, bounded
functions on $Q$ in $L_M(Q)^{3\times 3}_{\rm{sym}}$.  Then    $L_{M^\ast}(Q)^{3\times 3}_{\rm{sym}} = (E_M(Q)^{3\times 3}_{\rm{sym}})^*$ (cf. \cite{AW})
and we  observe that  $E_M\subseteq\mathcal{L}_M\subseteq L_M$.
The functional
    $$\varrho(\bbK)=\int_Q M(x,\bbK(x))dxdt$$
is a modular in the space of measurable functions
$\bbK:Q\to\R^{3\times 3}_{\rm{sym}}.$
\begin{definition}
We say that  sequence
$\{\tens{z}^j\}_{j=1}^{\infty}$ converges modularly to $\tens{z}$ in
$L_M(Q)^{3\times 3}_{\rm{sym}}$, which is denoted by  $\tens{z}^j\stackrel{M}{\to}\tens{z}$, if there exists
$\lambda>0$ such that
    $$\int_Q M\left(x,\frac{\tens{z}^j-\tens{z}}{\lambda}\right)dxdt\to0\quad\mbox{as }j\to\infty.$$
\end{definition}

\begin{definition}\label{delta2}
We say
that an $N$-function $M$ satisfies \emph{$\Delta_2$--condition} if for some
nonnegative, integrable on $\Omega$ function $g_M$ and a constant
$C_M>0$
	\begin{equation}
	M(x,2\bbK)\le C_M M(x,\bbK) + g_M (x)\quad\mbox{for all }\bbK\in\R^{d\times
	d}_{\rm{sym}} \mbox{ and a.a. } x\in\Omega.
	\end{equation}
\end{definition}
This condition is crucial for the structure of $L_M(Q)^{3\times 3}_{\rm{sym}}$  space.  It ensures that this space  is separable, reflexive and that   $C^\infty$  functions are dense.  What is more, if $\Delta_2$--condition 
holds, then ${E}_M(Q)^{3\times 3}_{\rm{sym}}=L_M(Q)^{3\times 3}_{\rm{sym}}$. Otherwise the considered space losses the above facilitating properties. 

\medskip

Below we recall several useful lemmas which are used within the  proof of existence of weak solutions.
Their proofs the reader can find e.g. in \cite{gwiazda-swierczewska,AW2013}.
\begin{lemma}\label{modular-conv}
Let $\tz^j:Q\to\R^{3\times 3}$ with  $ j = 1,\dots,\infty $ be  a measurable sequence. Then $\tz^j\modular{M}
\tz$ in $L_M(Q)^{3\times 3}$ modularly if and only if $\tz^j\to \tz$ in measure and
there exists some $\lambda>0$ such that the sequence
$\{M(\cdot,\lambda \tz^j)\}_{j=1}^{\infty}$ is uniformly integrable, i.e.,
$$\lim\limits_{R\to\infty}\left(\sup\limits_{j\in\N}\int_{\{(t,x):|M(x,\lambda\tz^j)|\ge
R\}}M(x,\lambda\tz^j)dxdt\right)=0.$$
\end{lemma}
	\begin{lemma}\label{uni-int}
	Let $M$ be an $N$--function and for all $j\in\N$ let $\int_Q
	M(x,\tz^j)dxdt\leq c$.
	Then the sequence $\{\tz^j\}_{j=1}^{\infty}$ is
	uniformly integrable.
	\end{lemma}
\subsection{Div--Curl lemma}
In further consideration we  use so called Div-Curl Lemma as given in \cite{tartar,FMR}. We denote for $\vec{a}= (a_0, a_1, a_2, a_3)$ 
\begin{equation} \label{divergence}
{\rm Div}_{t,x} \vec a  := \partial_t a_0 + \sum_{i=1}^3  \partial_{x_i}a_i
\quad
\mbox{ and } 
\quad
{\rm Curl}_{t,x} \vec a  := \nabla_{t,x} \vec a  -(\nabla_{t,x} \vec a)^T .
\end{equation}
\begin{lemma} 
Let $Q = (0,T) \times \Omega \subset  \R^4$ be a bounded set. Let $p,q,l,s  \in (1, \infty)$ be such that $ \frac{1}{p} +  \frac{1}{q}= \frac{1}{l}$ and vector fields $\vec a^n , \vec b^n$   satisfies 
\begin{equation*}
\vec a^n \rightharpoonup \vec a  \quad \mbox{weakly in } \quad L^p(Q)^4   \quad \mbox{and} \quad \vec b^n \rightharpoonup \vec b
\quad \mbox{weakly in } \quad L^q(Q)^4,
\end{equation*}
and $ Div_{t,x} \vec a$ and  $ Curl _{t,x} \vec b$ are precompact  in  $W^{-1,s} (Q)$  and $W^{-1,s} (Q) ^{4 \times 4}$ respectively. 
Then
\begin{equation*}
\vec a^n \cdot \vec b^n  \rightharpoonup \vec a \cdot \vec b  \quad \mbox{weakly in } \quad L^l(Q),
\end{equation*}
where $\cdot$ stands for scalar product in $ \R^4$.  
\label{divcurl}
\end{lemma}

\section{Main  result} \label{defsol}
We start with a definition of a weak solution of the problem \eqref{PP}.

    \begin{definition}\label{weak-for} Let $\vr_0$ satisfy \eqref{as-vr0}, $\uu_0\in H(\Omega)^3$, $ \theta_0 $ satisfy \eqref{as-th0} and $\bff \in L^{p'}(0,T;L^{p'}(\Omega)^3)$. Let  $\bS$ satisfy
	conditions {\bf{S1-S3}} with an $N$--function $M$ such that
	\begin{equation*}
		M(x,\txi)\ge \cp|\txi|^{p} - \widetilde{C}
		\quad \mbox{ with } \cp>0, \widetilde{C} \geq0 \mbox{ and }
		p\geq\frac{11}{5}   
		\end{equation*}
  for a.a. $x\in \Omega$ and all $\txi \in \R^{3 \times 3}_{\rm sym}$ and let $\vec q$ satisfies \eqref{zal_q} and \eqref{ass-q} with 
		$
		  \beta > - \min \left\{  \frac{2}{3}, \frac{3p -5 }{3p -3} \right \}.
		$

We call  ($\vr$,  $\uu$, $\theta$)  a \emph{weak solution} to (\ref{PP}) if
           \begin{equation*}
			\begin{split}
            & 0< \vr_* \leq \vr(t,x) \leq \vr^* \quad \mbox{for a.a. }(t,x)\in Q,\\
			& \vr \in C([0,T];L^q(\Omega))\quad\mbox{for arbitrary } q\in[1,\infty),\\
			& \partial_t\vr \in L^{5p/3}(0,T;(W^{1,5p/(5p-3)}(\Omega))^*),\\
			& \uu \in L^{\infty}(0,T;H(\Omega)^3)\cap L^p(0,T;\Vpd) \cap N^{1/2,2}(0,T; H(\Omega)^3),\\
            &{ \Du\in L_M(Q)^{3\times 3}_{\rm sym}}\quad\mbox{and}\quad(\vr\uu,\vec\psi)\in C([0,T])\mbox{ for all }\vec\psi\in H(\Omega)^3,\\
	& \theta \in L^{\infty}(0,T;L^1(\Omega)) \quad \mbox{ and }\quad \theta \geq \theta_*> 0  \mbox{ for } a.a \ (t,x) \in Q,  \\
	& \theta ^{\frac{\beta - \lambda +1 }{2}} \in L^{2}(0,T;W^{1,2}(\Omega)) \mbox{ for all  } \lambda \in (0,1),\\
	& \partial_t (\vr \theta) \in L^1(0,T; (W^{1,q}(\Omega))^*) \mbox{ with  $q$  sufficiently large }
			\end{split}
			\end{equation*}
   and  the following identities are satisfied: for continuity equation
        \begin{equation}\label{26}
        \int_0^T\left\langle \partial_t\vr,z\right\rangle - (\vr\uu,\Grad z) dt = 0
       \end{equation}
    for all $z\in L^r(0,T;W^{1,r}(\Omega))$ with $r=5p/(5p-3)$, i.e. 
    \begin{equation*}
    \int_{s_1}^{s_2}\int_\Omega\vr\partial_t z + (\vr\uu)\cdot\Grad z dxdt
     = \int_\Omega\vr z(s_2)-\vr z (s_1) dx
    \end{equation*}
    for all $z$ smooth and $s_1,\,s_2\in[0,T],\,s_1<s_2$; for the momentum equation 
        \begin{equation*}
        \begin{split}
        &-\int_{0}^T\int_{\Omega} \vr\uu\cdot\partial_t \bvarphi
        - \vr \uu\otimes\uu :\Grad\bvarphi
         + \Sdu : \grs \bvarphi dxdt  = \int_{0}^T\int_{\Omega} \vr\bff\cdot\bvarphi dxdt
        + \int_{\Omega}\vr_0\uu_0\cdot\bvarphi(0) dx
        \end{split}
        \end{equation*}
 for all $\bvarphi\in\D ((-\infty,T);\V)$; and for energy equation 
	 \begin{equation*}
\int_{0}^T\langle  \partial_t  (\vr\theta) ,  h \rangle
        - ( \vr \theta  \uu , \nabla h )  + (\qu  ,\nabla h ) dt  =  \int_{0}^T (\Sdu ,\Du h) dt  
        \end{equation*}
 for all $h \in L^\infty(0,T; W^{1,q}(\Omega))$ with $q$  sufficiently large.
Moreover, initial conditions are achieved in the following way
    	\begin{equation}\label{incond}
	\begin{split}
	& \lim\limits_{t\to 0^{+}}\|\vr(t)-\vr_0\|_{L^q(\Omega)} + \|\uu(t)-\uu_0\|^2_{L^2(\Omega)}=0\quad \mbox{ for\,arbitrary\,}q\in[1,\infty),\\
	& \lim\limits_{t\to 0^{+}} (\vr \theta(t), h) = (\vr_0 \theta_0,h) \quad \mbox{ for all } h\in L^\infty(\Omega).
	\end{split}
	\end{equation}
    \end{definition}
    
    \begin{theorem}\label{main}
	Let $M$ be an $N$--function satisfying for some $\cp>0,$ $\widetilde{C} \geq0$  and 
	for a.a. $x\in \Omega$ and all $\txi \in \R^{3 \times 3}_{\rm sym}$
		\begin{equation}\label{115}
		M(x,\txi)\ge \cp|\txi|^{p} - \widetilde{C}
		\quad 
\mbox{ with } 
		p\geq\frac{11}{5}.
		\end{equation}
Let us
	 assume that  the complementary function 
	\begin{equation}\label{cd2}
	M^* \mbox{ satisfies a }\Delta_2\mbox{-condition} \quad \mbox{ and }\quad \lim\limits_{|\tens{\xi}| \to \infty}\sup\limits_{x\in \Omega} \frac{M^*(x,\tens{\xi})}{|\tens{\xi}|} = \infty .
	\end{equation}
	Moreover, let  $\bS$ satisfy
	conditions {\bf{S1-S3}} and $\vec{q}$ satisfy \eqref{zal_q}, \eqref{ass-q} with 
	$ \beta > - \min \left\{  \frac{2}{3}, \frac{3p -5 }{3p -3} \right \}.$ Let $\uu_0\in H(\Omega)^3$, $\vr_0\in L^{\infty}(\Omega)$ with $0<\vr_*\leq\vr_0(x)\leq\vr^*<+\infty$ for a.a. $x\in\Omega$, $ \theta \in L^1(\Omega)$, $0< \theta_* \leq \theta _0$ for a.a. $  x \in \Omega$
	and 	$\bff\in L^{p'}(0,T;L^{p'}(\Omega)^3)$.
	Then there exists a weak solution to (\ref{PP}).
	\end{theorem}

In the following paper we consider the flow in the domain of space dimension $d=3$, just for the brevity of the paper. The existence result can be extended to the case of arbitrary $d\geq2$ and $p\geq\frac{3d+2}{d+2}$.


\section{Proof of the Theorem \ref{main} }\label{proof}
To prove the  Theorem \ref{main} we  proceed with  $n$--approximation of  problem \eqref{PP}.  First in order to show the existence of such  $n$--approximation we need to introduce additional two level approximation. The next step is to provide uniform estimates for $n$--approximate problem  which allow us to pass to the limit with $n \rightarrow  \infty$ and show weak sequential stability.

\subsection{Existence of the $n$--approximate problem}\label{naprox}

Let $\{\oo^i\}_{i=1}^{\infty}$ be a orthonormal  basis of $\Vpd$  such that $ \{ \oo_i \}_{i=1}^\infty \subset \Vdpd $ and elements of the basis are 
constructed  with the help of eigenfunctions of the problem
	$$(\!(\oo_i,\bvarphi)\!)_s=\lambda_i(\oo_i,\bvarphi)\quad\mbox{for
		all}\,\,\bvarphi\in V_s,$$
where
	\begin{equation}\label{Vsdef}
	V_s\equiv \mbox{ the closure of }{\mathcal V}
	\mbox{ w.r.t. the }W^{s,2}(\Omega)\mbox{-norm for } s >3,
	\end{equation}	
where  $(\!(\cdot,\cdot)\!)_s$ denotes the scalar product in $V_s$. Then the Sobolev embedding theorem provides
	\begin{equation}\label{embedding}
	W^{s-1,2}(\Omega)\hookrightarrow L^{\infty}(\Omega).
	\end{equation}
Then  approximate solution is given by 
\begin{equation}\label{un}
\uu ^n:= \sum_{i=1}^n \alpha ^n_i(t) \oo^i \quad \mbox{ for } i=1,2,\dots,
\end{equation}
with $\alpha^n_i \in C([0,T])$ and  the  triple  $(\vr^n, \uu^n,\theta^n)$ satisfies 
	 \begin{equation}\label{26a}
        \int_0^T\left\langle \partial_t\vr^n,z\right\rangle - (\vr^n\uu^n,\Grad z) dt = 0
        \quad \mbox{ for all }z\in L^r(0,T;W^{1,r}(\Omega)) \mbox{ with }r=5p/(5p-3),
       \end{equation}
 where  $\vr^n(0, \cdot) = \vr_0$ and
       \begin{equation}\label{nmom}
      \left \langle  \partial_t (\vr^n \uu^n ), \oo_i \right \rangle - (\vr^n \uu^n \otimes \uu^n, \Grad \oo_i) + (\Sdun, \tens{D} \oo_i) = (\vrn \vec{f}^n , \oo_i)
       \end{equation}
for all $i= 1,\dots, n$ and a.a. $ t \in [0,T]$, where $\uu^n(0,\cdot) = P^n \uu_0$ ($P^n$ denotes the projection of $H(\Omega)^3$ onto linear hull of $\{ \vec\omega_i \}_{i=1}^n)$ and 
	$$P^n \uu_0 \to \uu_0 \quad  \mbox{ strongly in } L^2(\Omega)^3,$$
	$$ \vec{f}^n \to \vec{f} \quad \mbox{ strongly in } L^{p'}(0,T; L^{p'} (\Omega)^3).$$ 
Moreover 	$ \theta^n \in   L^{\infty}(0,T;L^2(\Omega))\cup L^{s}(0,T;W^{1,s}(\Omega))$ with $s = \min  \left\{  2, \frac{5 \beta+10}{\beta +5} \right\}$ and $\theta^n \geq \theta_*$ in $Q$,
       \begin{equation}
\label{tempweak}
        \int_0^T\left\langle \partial_t(\theta^n  \vrn) ,h\right\rangle - (\theta^n \vrn\uu^n,\Grad h) + ( \kappa_0 \Grad \theta^n, \Grad h ) dt =     \int_0^T(\Sdun, \Dun h)  dt 
       \end{equation}
for all $h \in  L^{\infty}(0,T;W^{1,q}(\Omega)) $ for large $q$ and where $\theta^n(0,\cdot) = \theta^n_0$ s.t.
	$$\theta^n_0 \to  \theta_0 \quad\mbox{ strongly in }L^1(\Omega).$$

\subsubsection{Proof of the existence of the n-approximate problem} \label{proof_n_approx}
The existence of a triple ($(\vr, \uu, \theta) = (\vr^n, \uu^n,\theta^n)$) \footnote{In section \ref{proof_n_approx} we omit the superscript $n$ to simplify the notation.}  given by \eqref{un}--\eqref{tempweak} can be proven  by two-steps approximation. To this end we adopt the proof from  \cite[Section 6]{FMR} where more details can be found. Here we present only main steps of the reasoning for the convenience of the reader. The proof is based on standard artificial viscosity  technique.

In order to define the new two-step approximation let  us take $ \{  w_j \}^{\infty} _{j=1}$  a smooth basis  of $W^{1,2}(\Omega)$ orthonormal in $L^2(\Omega)$  spanning  the space where we construct a $k$-approximation of $\theta$. We look for the triple  $ (\vrke , \uke , \tke )$  where $\uke$ and $\tke$ are defined by
\begin{equation*}
\uke := \sum _{i=1} ^ n \ake_i (t) \oo_i   \quad   \mbox{ and  } \quad \tke := \sum  _{i=1}  ^ k \bke_i (t)  w_i  
 \end{equation*}
	and  $ (\vrke , \uke , \tke )$ satisfies the following 
\begin{align}
& \partial_t \vrke + \Div (\vrke \uke)  - \epsilon \Delta  \vrke = 0 \mbox{ in } Q, \quad \Grad \vrke \cdot \vec{ n} = 0\quad \mbox{on}\quad [0,T] \times \partial \Omega ,  \label{65}\\
&\vr_* \leq \vrke \leq \vr^*\quad  \mbox{in} \quad Q , \quad \vrke(0,\cdot) = \vr_0 \quad\mbox{in} \quad \Omega \label{nrho},\\
 & ( \vrke \frac {d}{dt} \uke , \oo_i)  + (\vrke[\Grad \uke]\uke, \oo_i) + (\Ske, \tens{D} \oo_i)  - \epsilon (\Grad \vrke, [\Grad \uke] \oo_i)= (\vrke \vec{f}^n,\oo_i) \notag\\
& \mbox{ in } Q  \mbox{ and for all }  i =  1,2,\dots,n,\label{66}\\
&\uke(0,\cdot) = \uke_0 =  \sum _{i=1} ^ n \alpha^{k,\epsilon}_i (0)\oo_i  = P^n \uu _0  \mbox{ in } \Omega, \notag\\
  \end{align}
   \begin{align}
& ( \vrke \frac {d}{dt} \tke , w_i)  + (\vrke[\Grad \tke]\uke, w_i) +(\kappa_{k,\epsilon} \Grad \tke , \Grad w_i)  -\epsilon (\Grad \vrke \Grad \tke, w_i)=  (\Ske :  \tens{D} \uke, w_i)  \notag \\
& \mbox{ in } Q \quad  \mbox{for all }   i =  1,2,\dots,k, \label{67} \\
& \tke (0,\cdot) = \tke _0= \sum _{i=1} ^ k \nu^{k,\epsilon}_i(0) w_i   = P^k (\theta ^n_0)  \mbox{ in } \Omega,\label{68}
  \end{align}
where $\tke_{max} := \max_{(x,t)\in Q} \left ( \tke, \theta_*\right )$, $ \Ske := \tens{S} (\vrke, \tke_{max}, \Du^ {k,\epsilon} ) \quad  \mbox{and} \quad  \kappa_{k,\epsilon} := \kappa_0 (\vrke,\tke_{max})$. Here  $P^n$ denote the projection of $H(\Omega)^3$ onto linear  hull spanned by  $\{ \oo_i \}_{i=1}^{n}$ and $P^k$ analogously projection of $L^2$ onto $ {\rm span}\{  w_j \}^{k} _{j=1}$.

Multiplying  the $i$-th equation in \eqref{66} by $  \ake$, taking  sum over $i= 1,\dots,n $,  using  $L^2(\Omega)$ scalar product of  \eqref{65}  with $| \uke |^2 /2$  and integrating over $(0,t)$ we obtain  that
\begin{equation}\label{momke}
\| \uke (t)\|_{L^2(\Omega)}^2+ \| \sqrt{\vrke}\uke (t)\|_{L^2(\Omega)}^2 + \int_0^t  M^*(x ,\Ske) + M(x, \tens{D} \uke) d \tau \leq C( \vr_0,\uu_0, \vec{f}),
\end{equation}
 what combined assumptions on tensor $\bS$ and with the  Korn inequality gives 
\begin{equation*}
\| \uke (t)\|_{L^{\infty}(0,T;L^2)} +   \| \uke (t)\|_{L^p(0,T;W^{1,p}(\Omega)^3)}  \leq C.
\end{equation*}
Furthermore, multiplying the  $i$-th equation of  \eqref{67} by $\bke_i $, taking sum over $i$, using  $L^2$ scalar product of  \eqref{65}  with $| \tke |^2 /2$ and integrating over $(0,t) $ leads to 
\begin{align*}
\| \tke (t)\|_{L^{\infty}(0,T;L^2)}+ \| \sqrt{\vrke}\tke (t)\|_{L^{\infty}(0,T;L^2)} + \| \sqrt{\kappa_{k,\epsilon} }\nabla \tke\|_ {L^{2}(0,T;L^2)}\\ \leq 
C  \| \sqrt{\vr_0 }\theta_0 \|_{L^2(\Omega)}^2 + C \|\Ske : \tens{D}(\uke) \|_ {L^{2}(0,T;L^2)} \leq C(n),
\end{align*}
where the last inequality holds because of  \eqref{momke} and the fact that $\uke$ is a linear combination of $n$ first elements of the basis $\{\oo^i\}_{i=1}^{\infty}$ such that $\Grad \oo^i \in L^{\infty }(\Omega)$  (see  \eqref{embedding}), hence $\|\Ske\|_{L^{\infty}}, \|\Du_{k,\epsilon}\|_{L^{\infty}} \leq C(n)$.

Applying once again  previous reasoning  but multiplying \eqref{66} by $ \frac{d \ake_i}{dt}$ and \eqref{67} by $ \frac{d \bke_i}{dt}$ we conclude that 
\begin{equation}
\label{wspke}
\begin{split}
&\|\ake_i \| _{W^{1,2}(0,T) } \leq C(n)  \mbox{ for  } i= 1,\cdots, n, \\
& \|\bke_i \| _{W^{1,2}(0,T) }   \leq C(k ) \mbox{ for  } i= 1,\cdots, k. 
\end{split}
\end{equation}

Summarising estimates \eqref{momke}--\eqref{wspke} together with \eqref{nrho} we can pass to the limit with $\epsilon$ and  obtain
 \begin{equation*}
\label{ke}
\begin{split}
&\vrke \weakstar \vr^k \quad \mbox{weakly-(*) in } L^{\infty} (Q),\\
&\ake_i\rightharpoonup \alpha^k  \quad  \mbox{weakly in } W^{1,2} (0,T) \mbox{ and strongly in } C([0,T]) \mbox{ for  } i= 1,\cdots, n, \\
&\bke_i  \rightharpoonup  \nu^k  \quad \mbox{weakly in } W^{1,2} (0,T) \mbox{ and strongly in } C([0,T])\mbox{ for  } i= 1,\cdots, k,\\
&\uke \rightarrow \uu ^k  \quad \mbox{strongly in } L^{2p} (0,T;W^{1,2p}_n(\Omega)),\\
&\vrke \rightarrow \vr^k  \quad  \mbox{strongly in } L^{2} (0,T;L^2(\Omega)) \mbox{ and  a.e. in } Q. 
\end{split} 
  \end{equation*}
Here $W^{1,2p}_n(\Omega)$ stands for $P^n(W^{1,2p})$. This set of convergence allows us to take the limit in \eqref{65}-\eqref{68} as $\epsilon \rightarrow 0$ and obtain that the  limit triple $(\vrk , \uk, \tk)$ solves
 \begin{align}
&\partial_t \vrk + \Div (\vrk \uk) = 0 \quad \mbox{ in }Q,  \quad \vr_* \leq \vrk \leq \vr^*\quad  \mbox{ in } Q , \quad \vrke(0,\cdot) = \vr_0 \quad\mbox{ in } \Omega ,\label{622}\\ 
 & ( \vrk \frac {d}{dt} \uk, \oo_i)  + (\vrk[\Grad \uk]\uk, \oo_i) + (\Sk, \tens{D} \oo_i) = 0 \quad\mbox{ in }  Q 
 \mbox{ for all } i =  1,2,\dots,n, \label{625}\\
& \uk(0,\cdot)  = P^n \uu _0  \quad\mbox{ in } \Omega, \notag\\
& ( \vrk \frac {d}{dt} \tk , w_i)  + (\vrk[\Grad \tk]\uk, w_i) +(\kappa_{k} \Grad \tke , \Grad w_i) =  (\Sk :\tens{D} \uk, w_i)   \quad\mbox{ in } Q  \mbox{ for all }   i =  1,2,\dots,k,  \notag\\
& \tk (0,\cdot) =  P^k  (\theta ^n_0) \quad\mbox{ in } \Omega . \label{628}
  \end{align}

The next step is to pass with $k \rightarrow \infty$ proceeding as in \cite{FMR}. Repeating procedures \eqref{momke}-\eqref{wspke} we obtain the following estimates 
\begin{align}
\| \uk\|_{L^{\infty}(0,T;L^2)} +   \| \uk\|_{L^p(0,T;W^{1,p})}  \leq C,\notag\\
\| \tk \|_{L^{\infty}(0,T;L^2)}+ \| \sqrt{\kappa_{k} }\nabla \tk\|_ {L^{2}(0,T;L^2)}\leq C , \notag \\  
\|\alpha^k_i \| _{W^{1,2}(0,T) } \leq C(n) \mbox { for } 1= 1,\dots, n. \notag 
\end{align}
Which imply that
\begin{equation}\label{uk}
\uk \rightarrow \uu  \mbox{ strongly in } L^{2p} (0,T;W^{1,2p}_n(\Omega)).
\end{equation}
Using the theory of renormalised solutions as in \cite{lions} we conclude that  for all $ p\in [1, \infty)$ holds
\begin{equation*}
\vrk \rightarrow \vr \mbox{ strongly in } C (0,T;L^p(\Omega)) \mbox{ and  a.e. in } Q. 
\end{equation*}
In addition, using procedure from  \cite{BMR2009} for selected subsequence we can  obtain that 
\begin{equation*}
\tk \rightarrow \theta \mbox{ strongly in } L^2(Q)  \mbox{ and  a.e. in } Q
\end{equation*}
and 
\begin{equation}\label{qk}
\kappa_{k} \nabla \tk  \weak  \kappa_0 \nabla \theta \mbox{ weakly in  } L^{\gamma }(0,T;W^{1,\gamma}(\Omega)) \mbox { for } \gamma = \min \left\{ 2, 1+ \frac{5}{(3 \beta +5 ) }\right\}. 
\end{equation}
  Summarising \eqref{uk}-\eqref{qk}  we pass to  the limit in the system \eqref{622}-\eqref{628} obtaining \eqref{26a}-\eqref{tempweak}. 
In addition from the minimum principle (see \cite[Section 6.4]{FMR}) we get 
\begin{equation}
\label{tng}
0< \theta_* \leq \theta^n \quad   a.a. \in Q. 
\end{equation}

\subsection{Uniform estimates}
Now we concentrate on passing with $ n \rightarrow \infty$. In first several steps we adopt reasoning from \cite{AW2013} and \cite{FMR}.  By standard method of characteristics for the transport equation (more details the reader can find in \cite{lions1,AW2013})
 we obtain
    \begin{equation}\label{bound-ro}
    0<\vr_*\leq\vr^n (t,x)\leq\vr^*<+\infty\quad{\rm for\,a.a.\,} (t,x)\in Q.
    \end{equation}
Multiplication  \eqref{nmom} by $\alpha_i^n$, taking sum up over $i= 1,\dots, n$ and use of  (\ref{26a}) leads to
	\begin{equation}\label{nn}
	\frac{1}{2}\frac{d}{dt}\int_\Omega\vr^n |\uu^n|^2dx + (\Sdun,\Dun) = (\vr^n\bff^n,\uu^n).
	\end{equation}
Using the H\"{o}lder, the Korn and the Young inequalities, assumption \eqref{115}  and   inequality~\eqref{bound-ro}
we are able to estimate the right-hand side of \eqref{nn} in the following way
    \begin{equation}\label{rhsnn}
    |(\vr^n \bff^n,\uu^n)|\leq 
 C(\Omega,c_c,\cp,\vr^*,p)  \left( \|\bff^n\|^{p'}_{L^{p'}(\Omega)}
    +\int_\Omega M(x,\Dun)dx \right).
    \end{equation}
Integrating (\ref{nn}) over the time interval $(0,s_0)$, using estimates
(\ref{rhsnn}), (\ref{bound-ro}),  the coercivity conditions {\bf {S2}} on $\bS$,
uniform continuity of $P^n$  w.r.t. $n$ and strong convergence $\bff^n \to \bff$ in $L^{p'}(0,T;L^{p'}(\Omega)^3)$
 we obtain
	\begin{equation}\label{osz}
	\begin{split}
	\int_\Omega \frac{1}{2}\vr^n(s_0)|\uu^n(s_0)|^2 dx&+\int_0^{s_0}\int_\Omega
	\frac{c_c}{2} M(x,\Dun)+ c_c M^\ast(x,\bS(t,x,\vr^n,\Dun))dxdt\\
	& \leq C (\Omega,c_c,\cp,\vr^*,p,\|\bff\|_{L^{p'}(0,T;L^{p'}(\Omega))})
	+\frac{1}{2}\vr^*\|\uu_0\|_{L^2(\Omega)}^2,
	\end{split}
	\end{equation}
where $C$ is a nonnegative constant independent of $n$ and dependent on the given data.

By \eqref{osz}, the condition (\ref{115}) provides that
$\{\Dun\}_{n=1}^{\infty}$ is uniformly bounded in the space $L^p(Q)^{3\times 3}$, i.e. 
	\begin{equation}\label{osz-gradSymLp}
	\int_0^T\|\Dun\|^p_{L^p(\Omega)}dt\leq C.
	\end{equation}
From the Korn inequality it is straightforward to show that 
	\begin{equation}\label{osz-gradLp}
	\int_0^T\|\Grad \uu^n\|^p_{L^p(\Omega)}dt\leq C.
	\end{equation}
Using (\ref{osz}) 
one can  deduce that
	\begin{equation}\label{osz-SDL1}
	\|\Sdun:\Dun \|_{L^1(Q)}\leq C,
	\end{equation}
	\begin{equation}\label{osz-SL1}
	\|\Sdun\|_{L^1(Q)}\leq C.
	\end{equation}
What is more, 
 the sequence $\{\Sdun\}_{n=1}^{\infty}$ is
uniformly bounded in Orlicz class ${\mathcal L}_{M^*}(Q)^{3\times 3}$.

Furthermore  (\ref{osz}) and (\ref{bound-ro}) provide
	\begin{equation}\label{313}
	\sup\limits_{t\in[0,T]} \|\uu^n(t)\|^2_{L^2(\Omega)}\leq C  \quad \mbox{ and } \quad  \sup\limits_{t\in[0,T]} \|\vr^n(t)|\uu^n(t)|^2\|_{L^1(\Omega)}\leq C ,
	\end{equation}
where $C$ is a positive constant dependent on the size of data,
but independent of $n$.
Since the sequence $\{\uu^n\}_{n=1}^{\infty}$ is uniformly bounded in $L^p(0,T;\Vpd)$ the Gag\-lia\-rdo-Nirenberg-Sobolev inequality gives us uniform bound
 in $L^p(0,T; L^{3p/(3-p)})$. Interpolation
  (see e.g. \cite[Proposition~1.41]{roubicek}) between spaces $L^\infty(0,T;L^2)$
and $L^p(0,T; L^{3p/(3-p)})$ provides
    \begin{equation}\label{osz-u5p3}
    \int_0^T\|\uu^n\|^{r}_{L^{r}(\Omega)}dt\leq C\quad \mbox{for }1\leq r\leq 5p/3
    \end{equation}
   (the above particular argument deals with the case $p<3$, the case $p \geq 3$ can be treated
 easier, e.g. with the Poincar\'e or the Morrey inequalities).
Therefore from (\ref{bound-ro}) and (\ref{osz-u5p3}) we infer also
    \begin{equation}\label{5p3est}
    \int_0^T\|\vr^n\uu^n\|^{5p/3}_{L^{5p/3}(\Omega)}dt\leq C.
    \end{equation}
Use of  (\ref{bound-ro}), (\ref{osz-gradLp}) and \eqref{osz-u5p3} combined with the H\"{o}lder inequality leads to 
\begin{equation*}
\int_0^T|(\vr^n\uu^n\otimes\uu^n,\Grad\uu^n)|dt\leq C  \quad\Longleftrightarrow\quad p\geq \frac{11}{5}. 
\end{equation*}
One can notice that here restriction for the exponent $p$ stated in \eqref{115} is given. 
By  (\ref{5p3est}) and  (\ref{bound-ro})   we obtain from \eqref{26a} that
\begin{equation}\label{317_1}
\int_0^T\|\partial_{t}\vr^n\|_{(W^{1,5p/(5p-3)})^*}^{5p/3}dt\leq C .
\end{equation}
\begin{lemma}
\label{lem-nik}
The sequence $\{ \uu^n \}_{n=1}^\infty$ is uniformly bounded w.r.t. $n$ in  Nikolskii space  $N^{1/2, 2}(0,T; H(\Omega)^3),$ namely 
\begin{equation}
\label{nik}
\sup\limits_{0<\delta<T} \delta^{-\frac{1}{2}} \left( \int_0^{T-\delta} \|\uu^n(s+\delta)-\uu^n(s)\|_{L^2(\Omega)} ds \right)^{\frac{1}{2}}  < C. 
	\end{equation}
	\end{lemma}
The proof of the above lemma can be found in   \cite[Section 3.1]{AW2013} (equations (55)-(62) therein). 
It is based on  reasoning from \cite[Chapter3. Lemma 1.2]{kazhikhov}  with modifications concerning a change of $L^2$ to $L^p$ structure and due to the nonlinear term controlled by the non-standard conditions \eqref{coercivity}. We notice that the presence of temperature does not influence this proof. The above lemma leads to the conclusion that $ \uu \in N^{1/2, 2}(0,T; H(\Omega)^3)$.

\medskip 


Finally we need to provide estimates concerning energy equation and the temperature. First we notice that taking $h=1$ in \eqref{tempweak}, by the Fenchel-Young inequality, \eqref{osz} and \eqref{bound-ro}   one gets
\begin{equation}
\label{tempsup}
\sup_{t \in [0,T]} || \vr^n \tn ||_{L^1(\Omega)} \leq C \quad \mbox{and}  \quad \sup_{t \in [0,T]} ||  \tn ||_{L^1(\Omega)} \leq C.
\end{equation}

Now, let us take $h = -(\tn) ^{-\lambda}$ with $\lambda \in (0,1)$ in \eqref{tempweak}. As $\tn \geq \theta_* $  (see \eqref{tng}) we have that $\|-(\tn) ^{-\lambda} \|_{\infty, Q} \leq C $ and from this substitution we obtain 
\begin{equation}
\label{temp1}
\int_0 ^T \| (\tn)^{\frac{\beta-\lambda-1}{2}} \Grad \tn \|^2_{L^2(\Omega)} dt  =  C_1 \int_0 ^T \| \Grad[ (\tn)^{\frac{\beta-\lambda+1}{2}}]\|^2_{L^2(\Omega)} dt \leq C_2 ,
\end{equation}
which provides
\begin{equation}
\label{tempw12}
\begin{split}
&\int_0 ^T \| (\tn)^{\frac{\beta-\lambda+1}{2}}  \|^2_{W^{1,2}(\Omega)} dt  \leq    \int_0 ^T \|   (\tn)^{\frac{\beta-\lambda+1}{2}}\|^2_{L^2(\Omega)}  dt + \int_0 ^T \| \Grad  (\tn)^{\frac{\beta-\lambda+1}{2}}\|^2_{L^2(\Omega)} dt	\\			 \underset{\mathrm{\eqref{temp1}}}{\leq}  
& \int_0 ^T \|   (\tn)^{\frac{\beta-\lambda+1}{2}}\|^2_{L^2(\Omega)} dt + C_1  \leq C_2 \Big[  \int_0 ^T \|     ( \tn)^{\frac{\beta-\lambda+1}{2}}\|^2_{L^1(\Omega)}  dt + \int_0 ^T \| \Grad  (\tn)^{\frac{\beta-\lambda+1}{2}}\|^2_{L^2(\Omega)} dt \Big] \\&+ C_1 \leq  
C_3 .
\end{split}
\end{equation}
From continuous embedding   of $ W^{1,2} $ in $ L^6$  we obtain 
$\int_0 ^T \| (\tn) \|^{\beta-\lambda+1}  _{3(\beta-\lambda+1)}    \leq   C$.
By interpolation with  \eqref{tempsup} we conclude that 
\begin{equation}
\label{temps}
\int_0 ^T \| (\tn) \|^{s}  _{L^s(\Omega)}   dt \leq   C   \quad \mbox{ for all  } s \in \left[1, \frac{5}{3}+\beta\right). 
\end{equation}
By the  assumption  made on heat flux \eqref{ass-q}  we have 
\begin{equation*}
| \kappa _0 \nabla_x \tn | \leq \kappa^*(\tn) ^{ \frac{\beta-\lambda-1}{2}} |\nabla  \tn| \  (\tn)^{\frac{\beta-\lambda+1}{2}}.
\end{equation*}
Let us notice that $ \| (\tn) ^{ \frac{\beta-\lambda-1}{2}} |\nabla  \tn|  \|_ {2} \leq C$ by \eqref{tempw12} and    $ (\tn) ^{ \frac{\beta-\lambda+1}{2} }  \in L^{\frac{2}{\beta-\lambda+1 }  ( \frac{5+3\beta}{3} - \varepsilon)  }$ by  \eqref{temps} with arbitrary small $\varepsilon >0$. Thus we see that
\begin{equation}
\label{tempkap}
\int_Q | \kappa _0 \nabla_x \tn |^m dx dt \leq C  \quad \mbox{for all } \quad m \in \left[1, \frac{5+3\beta}{4+3\beta}\right).
\end{equation}
Now we are ready to estimate last term in energy equation. By the Sobolev embedding, the Riesz-Torin interpolation theorem, the  H\" older inequality and above considerations we obtain
\begin{equation}
\label{tempfin}
\begin{split}
&\int_0 ^T \| \tn  \vrn \un \|^{\gamma}  _{\gamma}    dt  \leq \vr^* \int_0 ^T \| \un \|^{\gamma}  _{\frac{3p}{3-p}}  \|\tn\|^{\gamma}_{\frac{3p\gamma}{(3+\gamma)p-3\gamma}} dt\\& \leq 
C_1 \int_0 ^T \| \un \|^{\gamma}  _{1,p}  \|\tn\|^{(1-\alpha) \gamma}_{1} \|\tn\|^{\alpha  \gamma}_{3(\beta-\lambda_1)}  dt 
\underset{\mathrm{\eqref{tempsup},\eqref{temps}}} {\leq}  
C_1 \int_0 ^T \| \un \|^{\gamma}  _{1,p}  \|\tn\|^{\alpha  \gamma}_{3(\beta-\lambda+1)}dt \\ &
{\leq} 
C_1 \left [  \int_0 ^T \| \tn \|^{\beta-\lambda+1}  _{3(\beta-\lambda+1)}     dt  \right ] ^ {\frac{\gamma \alpha }{\beta - \lambda +1 }} \left [ \int_0 ^T \| \un \|^{\frac{(\beta-\lambda+1) \gamma}{\beta-\lambda+1 - \alpha\gamma} }  _{1,p}  dt  \right ]^  {\frac{\beta - \lambda +1 - \alpha \gamma  }{\beta - \lambda +1 }}
\underset{\mathrm{\eqref{temp1} }}{\leq}  C_2     .   
\end{split}
\end{equation}
 The parameter  $\alpha$ is taken such that
 \begin{equation} 
\label{alpha}
\frac{(3+\gamma)p-3\gamma}{3p\gamma} = \frac {1-\alpha}{1} + \frac{\alpha}{3(\beta - \lambda +1)}.
\end{equation}
 Moreover, the last inequality in \eqref{tempfin} gives constrains  combining values of $ \beta,\  \alpha, \   \lambda,\  p$  and $ \gamma$, i.e.
$
\frac{(\beta - \lambda +1)\gamma}{\beta - \lambda +1- \alpha \gamma  } = p.
$
Using formula \eqref{alpha}  we claim that $\gamma >1 \Leftrightarrow \beta > -\frac{3p -5}{3p -3}$ which is the restriction we demand in Theorem~\ref{main}. To sum up we obtain that for $p<3$ and appropriate $\beta$  there exist $\gamma  >1$ such that
$$\|\vrn\un \tn \| _{ L^1 (0,T; L^{\gamma}(\Omega))}<C.$$ The above result holds also for $p \geq 3 $  because of  embedding properties of $W^{1,p}$. 

In the end we obtain from \eqref{tempweak} and estimates \eqref{osz-SDL1}, \eqref{tempkap} and the fact that $\vrn\un \tn \in L^1 (0,T; L^{\gamma})$  that 
\begin{equation*}
\| \partial_t (\vrn \tn)\| _{ L^1 (0,T; (W^{1,s})^*)} < C  \quad \mbox{ for } s \mbox{ sufficiently large}.
\end{equation*}

The above considerations provide all necessary uniform estimates.

\subsection{Weak limits}
Uniform estimates obtained in previous section together with  the  Banach-Alaoglu theorem  ensure existence of subsequences selected from $\{ \vrn \}_{n=1}^{\infty}$, $\{\un \}_{n=1}^{\infty} $, $ \{\tn \}_{n=1}^{\infty}  $ such that 
\begin{align}
&\vrn \rightharpoonup \vr  \weakly  L^q(Q)  \mbox{  for any } q \in [1,\infty ) \mbox{ and weakly-(*) in  }  L^{\infty}(Q) ,  \label{rhoubound}\\
&\partial_t \vrn \weak \partial_t \vr \weakly L^{5p/3}(0,T;W^{1,5p/(5p-3)})^*),  \label{rhobound}\\
&   \un \weak\uu   \weakly L^{p}(0,T;W^{1,p}_{0,{\rm div}}(\Omega)^3)\mbox{ and }  L^{5p/3} (Q)^3 \mbox{ and weakly-(*) in  }  L^{\infty}(0,T; H(\Omega)^3), \label{0p} \\
& \tn \weak \theta \weakly  L^q (Q)  \mbox{ for any } q \in [1, 5/3+\beta), \label{thetaweak}\\
& \tn  \geq \theta_* >0  \mbox{ for a.a. } (t,x) \in Q.  \label{thetaogr}
\end{align}
In addition,  there exist $\overline{\vr \uu} \in L^{5p/3} (Q)^3$ and $\overline{\theta^{\alpha} }\in L^2(0,T;W^{1,2}(\Omega))$ such that
\begin{align}
&\vrn \un \weak \overline{\vr \uu }   \weakly  L^{5p/3} (Q)^3 , \label{rhoweak}\\ 
& (\tn)^{\alpha} \weak \overline{(\theta)^{\alpha} } \weakly L^2(0,T; W^{1,2}(\Omega)) \mbox { for } \alpha \in ( 0, (\beta+1)/2).   \label{thetaalpha}
\end{align}

What is more, noticing that  $E_M$ and $E_{M^*}$ are separable spaces and $(E_{M})^* = L_{M^*}$,     $(E_{M^*})^* = L_M$ 
we obtain
\begin{align}
	\grs \uu^n \weakstar \grs \uu \quad \mbox{weakly-(*) in }L_M(Q)^{3 \times 3}_{\rm sym}, \label{wsdu}\\
	\bS(\cdot,  \vr^n,\theta^n, \Dun) \weakstar \bchi \quad\mbox{weakly-(*) in } L_{M^*}(Q)^{3 \times 3}_{\rm sym},\label{wssm}
\end{align}
where $\bchi \in \mathcal{L}_{M^*} (Q)^{3 \times 3}_{\rm sym} $. Applying Lemma~\ref{uni-int} we conclude the uniform integrability
of   $\{\bS(\cdot,\vr^n,\theta^n,\Dun)\}_{n=1}^{\infty}$.
Consequently there exists a tensor
 $\bchi\in L^1(Q)^{3\times 3}$ such that
	\begin{equation}\label{01}
	\bS(\cdot,\vr^n, \tn,\Dun)\weak \bchi\quad \mbox{ weakly in } L^1(Q)^{3\times 3}.
	\end{equation}

\subsection{Strong convergence}
In this section we will  prove  strong convergence of the triple $(\vrn , \un , \tn)$ using the Aubins-Lions lemma  and the Lemma~\ref{divcurl}. 

Let us start with the velocity field. Since \eqref{313}$_1$, \eqref{osz-gradLp} and \eqref{nik} hold, due to  \cite[Theorem~3]{simon} we have
	\begin{equation*}
	\un \to \vec{u} \quad \mbox{ strongly in } L^2(Q)
	\end{equation*}
and by \eqref{0p}
	\begin{equation}\label{ustrong}
	\un \to \vec{u} \quad \mbox{ strongly in } L^q(Q) \mbox{ with }q\in \left[ 1, {5p}/{3} \right) .
	\end{equation}
Using \eqref{bound-ro}, \eqref{rhobound}  and the Aubin-Lions argument (see \cite{simon})  we obtain
\begin{equation*}
\vrn \rightarrow \vr  \quad  \mbox{ strongly in } C([0,T];(W^{1,5p/(5p-3)}(\Omega))^*).
\end{equation*}
The concept of the renormalised solutions  (see \cite{lions} for details) leads to 
\begin{equation}\label{rholq}
\vrn \rightarrow \vr  \mbox{ strongly in } C([0,T];L^q(\Omega)) \mbox{ for all }  q \in [1, \infty)  \mbox{ and a.e. in Q}. 
\end{equation}
Moreover,  one can show also that 
\begin{equation*}
\lim_{t\rightarrow t_0^+}   \| \vr(t) - \vr(0) \|_{L^q(\Omega)} = 0 \mbox{ for all }  q \in [1, \infty) 
\end{equation*}
(what gives the first part of \eqref{incond}) and 
\begin{equation}
\label{rho0}
\lim_{t\rightarrow0^+}  ( \vr(t), |\uu (t)| ^2) = (\vr_0 , |\uu_0|^2)   ,
\end{equation}
for details see \cite[Section 3.4]{AW2013}. 
The above strong convergence  of the velocity field provides
\begin{equation}\label{pierrho}
\sqrt{\vrn}\un \rightarrow \sqrt{\vr}\uu \quad \mbox{ strongly  in }  L^2(\Omega)^3
\end{equation}
and
\begin{equation*}
\sqrt{\vrn}\un (t)  \rightarrow \sqrt{\vr}\uu (t)  \quad \mbox{ strongly  in }  L^2(\Omega)^3 \mbox{ for almost all $t<T$} .
\end{equation*}

Arguments  \eqref{ustrong} - \eqref{pierrho} together with  \eqref{bound-ro} leads to 
\begin{equation*}
\begin{split}
\vrn\un  \otimes \un   =(\sqrt{\vrn})(\sqrt{\vrn}\un)   \otimes \un  \weak {\vr}\uu \otimes \uu  \weakly  L^{\gamma '}(0,T;W^{-1,\gamma '}(\Omega))
\end{split}
\end{equation*}
for $\gamma$ sufficiently large, i.e. $ \frac{1}{q}+\frac{6}{5p} + \frac{1}{\gamma} <1$ with $ q \in [1,\infty)$. Density argument together with \eqref{0p} provides 
\begin{equation*}
\begin{split}
\vrn\un  \otimes \un   =(\sqrt{\vrn})(\sqrt{\vrn}\un)   \otimes \un  \weak {\vr}\uu \otimes \uu  
\weakly  L^{p '}(0,T;W^{ -1,p'}_{\rm div}(\Omega)).
\end{split}
\end{equation*}
Now we  show that 
\begin{equation}\label{vrun}
\vrn u^n_i   \weak  \vr u_i \weakly L^q(Q) \mbox{ for all } q\in [1,5p/6] \mbox{ and } i=1,2,3.
\end{equation}
We define $\vec a^n = (\vrn, \vrn u_1, \vrn u_2, \vrn u_3)$ and $ \vec b^n_i = (u^n_i,0,0)$ for $ i = 1,2,3$.  Using \eqref{rhoubound} and  \eqref{rhoweak}  we obtain convergence  $\vec a^n  \weak  ( \vr, \overline{\vr u_1},\overline{\vr u_2},\overline{\vr u_3}) $ in $L^q  (Q)^4 $ for $ q\in [1, 5p/3]$  and $\vec b_i \weak  (u_i,0,0)$ in $ L^{5p/3} (Q)^4$.  From the definition  \eqref{divergence}$_1$ and continuity equation we obtain 
\begin{equation*}
{\rm Div}_{t,x}  \vec a^n = \partial_t \vrn + \Div (\vrn \un) = 0. 
\end{equation*}
From definition \eqref{divergence}$_2$ we have
\begin{equation*}
{\rm Curl}_{t,x} \vec b^n_i =  \begin{pmatrix}
  0& \nabla u_i^n\\
- (\nabla u_i^n )^T& \mathcal O\\
 \end{pmatrix}
\end{equation*}
( $\mathcal O$ denotes here zero $ 3 \times 3$ matirx). Due to the fact that  $\un \in L^p(0,T;W^{1,p}(\Omega)^3) $ we obtain   $\nabla \un \in L^r  \hookrightarrow  \hookrightarrow   W^{-1,p} (Q) .$ Consequently  the assumptions of Lemma~\ref{divcurl}  are satisfied and from its statement  we obtain \eqref{vrun}.

Using the same tool as above the convergence of the  $  \{\vrn \tn \un\}_{n=1}^{\infty}$ can be shown. To this end we set 
\begin{equation*}
\vec a^n = (\vrn \tn, \vrn \tn u^n_1 +\kappa_0 \nabla \tn, \vrn \tn u^n_2 +\kappa_0 \nabla \tn, \vrn \tn u^n_3 +\kappa_0 \nabla \tn)
\end{equation*}
and $ \vec b^n= ((\tn)^{\alpha},0,0)$ with $\alpha \in (0,(\beta+1)/2)$. Inequalities  \eqref{tempkap}, \eqref{tempfin}, \eqref{thetaogr} and  \eqref{rholq} ensure that $\vec{a}^n$ converges weakly to $\vec a$ in $L^s(Q)$ for some $s>1$ close to 1 and $\vec b \weak (\overline{\theta^\alpha}, 0 ,0 ) $ in $ L^r(Q)$ for $r$ such that $\frac{1}{s}+\frac{1}{r} <1 $ (which is possible for  small $\alpha$ and due to condition \eqref{thetaweak}). In view of the energy equation it holds that 
\begin{equation*}
\begin{split}
{\rm Div}_{t,x} \vec a  & = \partial_t (\vrn \tn) + \Div(\vrn \tn \un +\kappa_0 \nabla \tn) 
= \Sdun : \Dun  \in L^1 (Q)\hookrightarrow \hookrightarrow W^{-1,\hat{r}}(Q),
\end{split}
\end{equation*}
where $ \hat{r} \in (1, 3/4)$.     On the other hand 
\begin{equation*}
{\rm Curl}_{t,x} \vec b^n =  \begin{pmatrix}
  0& \nabla (\tn)^{\alpha}\\
- (\nabla (\tn)^{\alpha} )^T& \mathcal O\\
 \end{pmatrix} \in L^2 (Q)^{4\times 4} \hookrightarrow \hookrightarrow (W^{-1,2}(Q)^{4\times 4} )^*.
\end{equation*}
The statement of Lemma~\ref{divcurl} provides that
\begin{equation*}
\vrn (\tn)^{\alpha +1 } \weak \vr \theta \overline{\theta^\alpha} \weakly L^{1+\eta}(Q) \mbox{ for some } \eta >0.  
\end{equation*}
The above combined with \eqref{temps} and \eqref{rholq} gives 
\begin{equation}\label{thet}
\vr (\tn)^{\alpha +1 } \weak \vr \theta \overline{\theta^\alpha} \weakly L^{1+\zeta}(Q) \mbox{ for some } \zeta >0.  
\end{equation}
Next step is  to show that  $ \overline{\theta^\alpha}  = \theta^{\alpha} $   a.e. in Q.   To do so we employ  Minty's trick. For $ y \in R^+$ and $ \alpha >0, y^{\alpha} $ is a increasing function which leads to 
\begin{equation*}
0 \leq \int_0^T  \big( \vr[ (\tn)^\alpha - h^{\alpha} ], \tn -h  \big) dt  \quad \mbox{ for all } h \in  L^{1+\eta}(Q).
\end{equation*}
Passing to the limit with $ n \to \infty $ and using  \eqref{thet} we obtain
\begin{equation*}
0 \leq \int_0^T  \big( \vr[ \overline{\theta^\alpha} - h^{\alpha} ], \theta -h  \big) dt \quad  \mbox{ for all } h \in  L^{1+\eta}(Q).
\end{equation*}
By setting  $ h = \theta - \lambda v $ for $\lambda >0,  v \in L^{1+\eta} (Q)$ and  $ h = \theta + \lambda v $  then passing to the limit with $ \lambda \rightarrow 0 $ we conclude that  

\begin{equation*}
0 = \int_0^T  \big( \vr[ \overline{\theta^\alpha} - \theta ^{\alpha} ],  v   \big) dt \quad  \mbox{ for all } v \in  L^{1+\eta}(Q).
\end{equation*} 
Therefore as $ \vr > \vr_*$ we deduce
\begin{equation*}
\overline{\theta^\alpha}  = \theta^{\alpha}   \quad  \mbox{ a.e. in Q. }  
\end{equation*}

Then by \eqref{thet}, weak convergence in $L^{1+\alpha}(Q) $ of $ \{ \vr ^{\frac{1}{1+\alpha } }\tn \}_{n=1}^{\infty}  $  to $ \vr ^{\frac{1}{1+\alpha } }\theta$
and convergence of $ \| \vr ^{\frac{1}{1+\alpha } }\tn  \| _{L^{1+\alpha}(Q)}$ to   $ \| \vr ^{\frac{1}{1+\alpha } }\theta  \| _{L^{1+\alpha}(Q)}$   are provided. Consequently 
\begin{equation*}
\vr ^{\frac{1}{1+\alpha }} \tn  \rightarrow \vr ^{\frac{1}{1+\alpha } }\theta  \quad \mbox{ strongly in }   L^{1+\alpha}(Q),
\end{equation*}
which combined with \eqref{thetaweak} and \eqref{rhoubound} leads to 
 \begin{equation}\label{thetaae}
 \tn  \rightarrow  \theta \quad  \mbox{ strongly in }   L^{q } \mbox{ for all }  q \in [1, 5/3 + \beta) \mbox{ and a.e. in Q}.
\end{equation}
The above strong convergence together  with  \eqref{thetaalpha} ensure that 
 \begin{equation}\label{tealpha}
 (\tn) ^{\alpha} \weak \theta^{\alpha}  \quad  \mbox{ weakly  in }   L^{2 }(0,T; W^{1,2}(\Omega))  \mbox{ for all }  \alpha  \in (0,  (\beta+1)/2 ). 
\end{equation}
Using   the same  arguments as in   \eqref{tempfin}, \eqref{ustrong}, \eqref{rholq} and \eqref{thetaae} we  conclude that 
\begin{equation*}
\vrn  \tn \un   \rightarrow  \vr \theta \uu  \quad  \mbox{ strongly in }   L^{1}(Q)^3. 
\end{equation*}
Next step is to establish convergence of $\qn$. According to \eqref{zal_q}  
\begin{equation}\label{remind}
\qn = \kappa_0(\vrn, \tn)\nabla_x \tn 
      = \frac{2}{\beta - \lambda +1}(\tn)^{   \frac{-\beta + \lambda +1}{2} } \kappa_0(\vrn, \tn) \nabla_x  (\tn) ^{\frac{\beta - \lambda +1}{2} }.
\end{equation}
Inequality   \eqref{temps} can be used to provide 
\begin{equation}  \label{thetafinal}
\int_0 ^T \| (\tn)^{   \frac{-\beta + \lambda +1}{2} }\kappa_0(\vrn, \tn)  \|^{2r}_{L^{2r}(\Omega)} dt  \leq \int_Q (\tn) ^{r(\beta + \lambda +1)} dxdt \leq C
\end{equation} 
for $r$ s.t. $r  (\beta + \lambda +1) =5/3 + \beta - \lambda$. Notice that $r>1 $  for    $ \lambda$ small enough. Then almost everywhere convergence of $\{\vrn\}_{n=1}^{\infty}, \{\tn\}_{n=1}^{\infty}$ showed in \eqref{rholq} and \eqref{thetaae} combined with Vitalli's convergence theorem  and   
\eqref{thetafinal} leads to
\begin{equation*}
(\tn)^{   \frac{-\beta + \lambda +1}{2} }\kappa_0(\vrn, \tn) \rightarrow \theta^{   \frac{-\beta + \lambda +1}{2} }\kappa_0(\vr, \theta)\quad \mbox{ strongly in } L^2 (0,T;L^2(\Omega)).
\end{equation*} 
Mowever, convergence proved in \eqref{tealpha}  gives us
\begin{equation*}
\nabla_x (\tn)^{   \frac{\beta - \lambda +1}{2} } \weak \nabla_x \theta^{   \frac{\beta - \lambda +1}{2} } \quad  \mbox{ weakly  in } L^2 (0,T;L^2(\Omega)^3),
\end{equation*} 
which applied to  \eqref{remind} and by \eqref{tempkap} leads to  
\begin{equation*}
\qn \weak \qu  \quad \mbox{ weakly  in } L^s (Q)^3  \mbox{ for all } s \in \left( 1, \frac{5+3\beta }{4+ 3 \beta} \right)  . 
\end{equation*} 

Arguments established in this section allows us to pass to the limit in system \eqref{26a}-\eqref{nmom}. 
It remains only to characterise the nonlinear term and to show convergence in the the RHS of the energy equations \eqref{tempweak}.


\subsection{Integration by parts}
\label{intbyparts}
Let us notice that classical integration by parts formula does  hold for our considered problem since 
Orlicz spaces are not reflexive and  smooth functions are not dense if $\Delta_2$--condition is not satisfied. In general  also  there is no equivalence between Bochner type space $L_M(0,T;L_M(\Omega))$ and $L_M(Q)$, which holds only in case $N$--function $M$ is of polynomial type (see \cite{Donaldson}). Therefore let us recall the following result from \cite{AW2013}:

\begin{lemma}
Let exponent $p $, function $\vec{f}$,  $N$--functions  $M$  and $M^*$     be as in  Theorem \ref{main}. We assume that 
\begin{equation*}
			\begin{split}
            & 0< \vr_* \leq \vr(t,x) \leq \vr^* \quad \mbox{for a.a. }(t,x)\in Q,\quad \vr \in C([0,T];L^q(\Omega))\quad\mbox{for arbitrary } q\in[1,\infty),\\
			& \partial_t\vr \in L^{5p/3}(0,T;(W^{1,5p/(5p-3)}(\Omega))^*),\\
			& \uu \in L^{\infty}(0,T;H(\Omega)^3)\cap L^p(0,T;\Vpd) \cap N^{1/2,2}(0,T; H(\Omega)^3),\\
            &{ \Du\in L_M(Q)^{3\times 3}_{\rm sym}}\quad\mbox{and}\quad(\vr\uu,\vec\psi)\in C([0,T])\mbox{ for all }\vec\psi\in H(\Omega)^3,\\
            & \bchi \in      L_{M^*}(Q)^{3\times 3}_{\rm sym}
          \end{split}
\end{equation*}
and the couple ($\vr$, $\uu $) is a weak solution of $\partial_t \vr + \Div (\vr \uu)=0$  (see Definition \ref{weak-for})  and  satisfies
 \begin{equation}\label{bypart1}
        -\int_{0}^T\int_{\Omega} \vr\uu\cdot\partial_t \bvarphi
        - \vr \uu\otimes\uu :\Grad\bvarphi
         + \bchi : \grs \bvarphi dxdt  = \int_{0}^T\int_{\Omega} \vr\bff\cdot\bvarphi dxdt
        \end{equation}
 for all $\bvarphi\in\D ((0,T);\V)$. Then for a.a. $s_0 $ and $s$ s.t. $ 0< s_0 \leq s$  it holds that 
     \begin{equation}\label{limh}
    \frac{1}{2}\int_\Omega \vr(s,x)|\uu(s,x)|^2 dx
      +\int_{s_0}^{s}\int_\Omega \bchi : \Du dxdt
      = \int_{s_0}^{s}\int_\Omega \vr\bff\cdot \uu dxdt
    + \frac{1}{2}\int_\Omega \vr(s_0,x)|\uu(s_0,x)|^2 dx.
    \end{equation}
\end{lemma}
The detail proof can be found in \cite[Section 3.3]{AW2013} and it is based on a proper choice of a test function in \eqref{bypart1} and goes via Steklov regularisation with respect to the time variable.

Let us remark  that now the assumptions of Lemma~ \ref{intbyparts} are satisfied due to the previous subsections, in particular
\eqref{limh} holds for sufficiently rich family of test functions by density arguments.


\subsection{Monotonicity method}
In this section we investigate the weak limit $\bchi $ and we  show that   $\bchi  = \Sdu$ a.e.  in $Q$. The proof is based on monotonicity method adopted to nonreflexive anisotropic Musielak-Orlicz spaces. Reasoning follows the one presented in \cite{GSW2,AW2013} with modifications which allow to deal with dependence of $\bS$ on density and temperature and  we recall it here  for the convenience  of the reader.

Using integration by parts formula (see Lemma \ref{intbyparts}) and letting $s_0 \rightarrow 0$ (see  \eqref{rho0}) we obtain that
    \begin{equation*}
    \frac{1}{2}\int_\Omega \vr(s,x)|\uu(s,x)|^2 dx
     +\int_{0}^{s}\int_\Omega \bchi : \Du dxdt
    = \int_{0}^{s}\int_\Omega \vr\bff\cdot \uu dxdt
    + \frac{1}{2}\int_\Omega \vr_0(x)|\uu_0(x)|^2 dx.
    \end{equation*}
After integrating equation \eqref{nn} over the interval $(0,s)$, passing with $n\to \infty $   and comparing the result  with the above one may  conclude that
    \begin{equation}\label{limsup}
    \limsup\limits_{n\to\infty}\int_0^s\int_{\Omega}\Sdun: \Dun dxdt \leq
   \int_0^s \int_{\Omega}\bchi : \Du dxdt.
    \end{equation}
Denoting by  $Q^s$  time-space cylinder  $(0,s)\times\Omega$ and using monotonicity of $\bS$ (see condition {\bf {S3}})  one obtains  that 
\begin{equation}\label{43}
    \int_{Q^s}(\bS(x,\vr^n,\theta^n, \vd)-\Sdun):(\vd-\Dun)dxdt \geq 0
    \end{equation}
holds  for all $\vd\in L^{\infty}(Q)^{3\times 3}$. Let us notice that  $\bS(x,\vr^n,\theta^n,\vd) \in
L^{\infty}(Q)^{3\times 3}$.  This statement can by proven  by contradiction. To do so suppose that $\bS(x,\vr^n,\theta^n,\vd)$
is unbounded. Since $M$ is nonnegative, by coercivity condition  (\ref{coercivity}),
it holds that 
    $$|\vd|\ge\frac{M^\ast(x,\bS(x,\vr^n,\theta^n,\vd))}{|\bS(x,\vr^n,\theta^n,\vd)|}.$$
Then the  right-hand side tends to infinity as $|\bS(x,\vr^n,\theta^n,\vd)|\to\infty$ by (\ref{limitinfty}), which
contradicts that $\vd \in L^{\infty}(Q)^{3\times 3}$. 

Employing continuity of 
$\bS$ w.r.t.  second and third argument and a.e. convergence of $\{ \vr^n\}_{n=1}^{\infty}$, $\{ \tn\}_{n=1}^{\infty}$   we obtain a.e. convergence  of the $\{\bS(x,\vr^n,\theta^n,\vd)\} _{n=1}^{\infty}$ to  $ \bS(x,\vr,\theta,\vd)$.  Since  $\{\bS(x,\vr^n,\theta^n,\vd)\} _{n=1}^{\infty} \subset L^{\infty}(Q^s)^{3 \times 3}$ we obtain uniform integrability of  $\{M^*(\bS(x,\vr^n,\theta^n,\vd))\}_{n=1}^\infty$. Therefore by  Lemma~\ref{modular-conv}   modular convergence of the sequence $\{\bS(x,\vr^n,\theta^n,\vd)\} _{n=1}^{\infty}$ is provided in $L_{M^*}(Q)^{3 \times 3}$.  As $M^*$ satisfies the $\Delta_2$--condition, then the 
modular and strong convergence in $L_{M^*}$ coincide (see \cite{Kufner}) and  
$\bS(x,\vr^n,\theta^n,\vd)\to\bS(x,\vr^n,\theta^n,\vd)$ strongly in $L_{M^*}$.
Therefore by \eqref{wsdu} we deduce
\begin{equation}\label{sro}
\lim\limits_{n\to\infty}\int_{Q^s}\bS(x,\vr^n,\theta^n,\vd) : \Dun dxdt = \int_{Q^s}\bS(x,\vr,\theta,\vd) : \Du dxdt.
\end{equation}
  Passing to the
limit with $n\to\infty$, by  \eqref{0p}, \eqref{wssm},  \eqref{limsup},  \eqref{sro}  we obtain from \eqref{43} that
    \begin{equation}\label{45}
    \int_{Q^s}\bchi : \Du dxdt \geq \int_{Q^s} \bchi :\vd dxdt +
    \int_{Q^s}\bS(x,\vr,\theta,\vd) : (\Du-\vd)dxdt
    \end{equation}
and consequently
    \begin{equation}\label{46}
    \int_{Q^s} (\bS(x,\vr,\theta^n,\vd)-\bchi):(\vd-\Du)dxdt \geq 0.
    \end{equation}
    Let us set 
     $$\vd=(\Du) \bbbone_{Q_{i}} + h\zd\bbbone_{Q_{j}},$$
    with
     $ Q_{k} = \{(t,x)\in Q^s : |\grs\uu(t,x)|\leq k \ \ \mathrm{a.e.\ in\ } Q^s\} $ and where
 $k>0$, $0<j<i$, $h>0$ and  $\zd\in L^{\infty}(Q)^{3\times 3}$ are arbitrary. 
As $\bS(x,\vr,\theta,\tens{0})=\tens{0}$, from  \eqref{46} we have 
    \begin{equation}\label{13}
    -\int_{Q^s \backslash Q_{i}} (\bS(x,\vr,\theta,\tens{0})-\bchi): \Du dxdt + h
    \int_{Q_{j}}(\bS(x,\vr,\theta,\Du+h\zd)-\bchi):\zd dxdt \geq 0
    \end{equation}
and obviously
    $$\int_{Q^s \backslash Q_{i}}\bchi:
    \Du dxdt=\int_Q(\bchi : \Du)\bbbone_{Q^s \backslash Q_{i}}dxdt.$$
 Since $\bchi \in \mathcal{L}_{M^*}(Q)^{3 \times 3}$, $ \Du \in  \mathcal{L}_{M}(Q)^{3 \times 3}$ (which is a consequence of convexity, nonnegativity  of $ M^*$,  $ M $  and of weak lower semi-contiunity and estimate \eqref{osz}) by the Fenchel-Young inequality  we obtain that 
    $\int_Q \bchi : \Du dxdt<\infty$ and consequently 
     $$(\bchi : \Du)\bbbone_{Q^s \backslash Q_{i}}
    \rightarrow 0 \quad\mathrm{\ a.e.\ in\ }Q \mbox { for } i \rightarrow\infty.$$
Hence by the Lebesgue dominated convergence theorem
    $$\lim\limits_{i\to\infty}\int_{Q^s \backslash Q_{i}}
    \bchi : \Du dxdt=0.$$
 Letting $i\rightarrow\infty$ in (\ref{13}) and dividing by $h$, we get
    $$\int_{Q_{j}}(\bS(x,\vr,\theta,\Du+h\zd)-\bchi) : \zd dxdt\geq 0.$$
Since $\Du+h\zd\rightarrow \Du$ a.e. in $Q_j$ when $h\rightarrow
0^{+}$ and as $\{ \bS(x,\vr,\theta,\Du+h\zd)\}_{h>0}\subset L^{\infty}(Q_j)^{3\times 3}$,
$|Q_j|<\infty$, by the Vitali lemma we conclude
    $$\bS(x,\vr,\theta,\Du+h\zd)\rightarrow\bS(x,\vr,\theta,\Du)\quad \mathrm{ in } \quad L^1(Q_j)^{3\times3}\mbox{\ as \ }  h \to 0^+$$
and
    $$\int_{Q_j}(\bS(x,\vr,\theta,\Du+h\zd)-\bchi) : \zd
dxdt\rightarrow\int_{Q_j}(\Sdu-\bchi) : \zd dxdt   \quad \mbox{ as }  h \to 0^+ . $$
Consequently,
    $$\int_{Q_{j}}(\Sdu-\bchi) : \zd dxdt \geq 0 \quad \mbox{ for all }\zd\in L^{\infty}(Q)^{3\times 3}.$$
 The choice
    of  $\zd$ s.t. $\zd=-\frac{\Sdu-\bchi}{|\Sdu-\bchi|}$ if  $
    \Sdu\neq\bchi$ and 
    $\zd = 0 $ if  $ \Sdu=\bchi$
yields
    $$\int_{Q_j}|\Sdu-\bchi|dxdt\leq0.$$
Hence $  \Sdu=\bchi $   a.e. in $Q_{j}$ and as   $j$ is arbitrary it holds  also a.e. in
$Q^s$  for almost all  $s$ such that $0<s<T$. 
Finally we conclude   that
\begin{equation}
\label{eq}
\bchi=\Sdu   \quad \mbox{ a.e. in } Q . 
\end{equation}

\subsection{Convergence of $\Sdun : \Dun$} \label{strongSdu}
The next  crucial  part of the proof  is to establish convergence of  the  sequence   $\{ \Sdun : \Dun\} _{n=1}^{\infty}$. The idea follows  \cite{gwiazda2015renormalized} (later used also in \cite{Klawe}) and is based on the concept of biting convergence (see \cite{ball}) and the   theory  of Young measures (for details see \cite{muller}).   Let us start with recalling definition of bitting limit   (see \cite{ball}).
\begin{definition} [Biting limit]
Let $\{a_n\}_{n=1}^{\infty}$ be a bounded  sequence in $L^1(Q)$. We say that $a \in L^1(Q)$ is a biting limit of subsequence of  $
\{a_n\}_{n=1}^{\infty}$ if there exists nonincreasing sequence     $\{E_k\}_{k=1}^{\infty},  E_k \subset Q $   satisfying $ \lim_{k \rightarrow \infty }|E_k| = 0$  such that $a^n$ converge weakly to $a$  in $L^1(Q \setminus  E_k)$. We denote biting convergence with $\biting$
\end{definition}
\begin{lemma}\label{bitlemma}
Let $\{a_n\}_{i=1}^{\infty} $ be a bounded sequence in $L^1(Q)$ and let $0\leq a_0 \in L^1(Q) $. If assumptions 
\begin{itemize}
\item[{\bf{A1)}}] $a_n\geq -a_0$ for all $ n  = 1 ,\dots, \infty$,
\item[{\bf{A2)}}] $a_n \biting a  $ as $ n \to \infty $,
\item[{\bf{A3)}}] $\limsup_{n\to \infty} \int_Q a_n   dx dt \leq \int_Q a dx dt$,
\end{itemize}
hold, then 
\begin{equation*}
a_n \weak a  \weakly   L^1(Q) \mbox{ as }  n\to \infty . 
\end{equation*}
\end{lemma}
For the proof see \cite{pedregal,gwiazda2015renormalized}.

We show now that for $ \{a_n\}_{i=1}^{\infty} = \{\Sdun : \Dun\}_{i=1}^{\infty} $ the assumptions of  lemma \ref{bitlemma} are fulfilled which lead to a weak convergence of $ a_n$ in the $L^1(Q)$ space.
Assumption {\bf{A1}}  is fulfilled  due the coercivity  condition 
\textbf{S3}, namely  $ a_n \geq 0$.  Then {\bf{A3}} is a straightforward consequence of   \eqref{limsup}. What is left, is the  {\bf{A2}} -- biting convergence of  $ a_n = \Sdun : \Dun$ to  $ a = \Sdu : \Du$.

Using the monotonicity of $\tens{S}$  (see \textbf{S3}) we can write down that
\begin{equation} \label{rhslambda}
0\leq  (\Sdun - \Sdurt):(\Dun-\Du).
\end{equation}
The right hand side of the above inequality is uniformly bounded in $L^1(Q)$ space,
which holds due to the H\"older inequality for Orlicz spaces and uniform estimates obtained in \eqref{osz}. In particular, the uniform boundedness of  $ \Sdurt $  in $ L_{M^*}$ is a consequence of the following reasoning:
By coercivity condition  \textbf{ S2}, the Fenchel-Young inequality and  convexity of the  $N$--function $M^*$ we can deduce that 
\begin{equation}
c M(x, \Du) + \frac{2c -d }{2} M^* (x, \Sdurt) \leq M(x , \frac{2}{d} \Du)
\end{equation}
with $d= \min \{c,1 \}$.  Then since $\Du \in L_{M}(Q)^{3\times 3}$, we obtain that $ \{ \Sdurt\} _{n=1}^{\infty} $ is uniformly bounded in $L_{M^*}(Q)^{3\times 3} $.  Since RHS of  \eqref{rhslambda} is uniformly bounded in $L^1(Q)$, there exists a Young measure  $\mu_{x,t}(\cdot,\cdot)$ (see \cite[Theorem 3.1]{muller}) satisfying 
\begin{equation}\label{3g}
\begin{split}
 & (\Sdun - \bS(x,\vr^n,\tn ,\Du)):(\Dun-\Du) \notag \\   
 &\to_b  \int_{\R^2\times\R^{3\times 3 }} ( \Ssa  -  \bS(x,l,s,\Du))):(\lambda - \Du)d\mu_{x,t}(s,l,\lambda) 
\end{split}
\end{equation}
as $n \rightarrow \infty$. 
Using \cite[Corrolary 3.4]{muller} and \eqref{rholq} together with  \eqref{thetaae} we have that in fact    $\mu_{x,t}(\cdot,\cdot,\cdot)$ can be written down in a form $ \delta_{\theta,\vr}(l,s) \otimes \nu_{x,t}(\lambda)$. This leads to 
\begin{equation}\label{lastlambda}
\begin{split}
 &\int_{\R^2\times\R^{3\times 3 }} ( \Ssa  - \bS(x,l,s,\Du)):(\lambda - \Du)d\mu_{x,t}(l,s,\lambda)\\&=
 \int_{\R^{3\times 3 }} ( \Ssl  - \Sdu):(\lambda - \Du) d \nu_{x,t}(\lambda) \\&=
 \int_{\R^{3\times 3 }} \Ssl: (\lambda - \Du) d \nu_{x,t}(\lambda)- 
  \int_{\R^{3\times 3 }} \Sdu: (\lambda - \Du) d \nu_{x,t}(\lambda).
\end{split}
\end{equation}
On the other hand 
\begin{align}
\label{4g}
\int_{\R^{3\times 3 }} \Sdu: (\lambda - \Du) d \nu_{x,t}(\lambda)= 
\Sdu : \left(  \int_{\R^{3\times 3 }} \lambda d \nu_{x,t}(\lambda) - \Du \right)  =0. 
\end{align}
Indeed, the above holds since   $\Sdu $ is independent of $\lambda$. Moreover   $
 \int_{\R^{3\times 3 }} \lambda d \nu_{x,t}(\lambda) = \Du$  for a.e. $(t,x)\in Q$ by \cite[Theorem 3.1]{muller} and $\Dun \weak \Du$ in  $L^1(Q)^{3 \times 3}$ (consequence of \eqref{wsdu}). Then the second therm of RHS of  \eqref{lastlambda} disappears and  \eqref{lastlambda}  becomes
\begin{equation}
\label{g1}
\begin{split}
 &\int_{\R^2\times\R^{3\times 3 }} ( \Ssa  - \bS(x,l,s ,\Du )):(\lambda - \Du)d\mu_{x,t}(l,s,\lambda)
 =\int_{\R^{3\times 3 }} \Ssl: (\lambda - \Du) d \nu_{x,t}(\lambda).
\end{split}
\end{equation}
Furthermore,  as  $\{a_n\}_{n=1}^{\infty}$  is uniformly bounded in $L^1 (Q)$ (by Fenchel-Young inequality and \eqref{osz})
we get that 
\begin{align} 
\Sdun:\Dun \to_b&  \int_{\R^2\times\R^{3\times 3 }} \Ssa: \lambda  d \mu_{x,t}(l,s,\lambda) \notag\\
&=\int_{\R^{3\times 3 }}\Ssl : \lambda d \nu_{x,t}(\lambda) .
\notag \end{align}
 Then as $ a_n \geq 0 $ for $ n= 1,\dots , \infty $, by  \cite[Corrolary 3.3]{muller}   and by  \eqref{limsup}, \eqref{eq},    we obtain  that 
\begin{align}
\label{2g}
\int_Q  \Sdu:\Du dxdt  &\geq  \liminf_{n \to \infty}  \int_{Q} \Sdun : \Dun  dx dt \\
&\geq \int_Q \int_{\R^{3\times 3 }} \Ssl: \lambda  d \nu_{x,t}(\lambda) dxdt. \notag 
\end{align}
However  $ \Sdu = \int_{\R^{3\times 3 }} \Ssa  d \nu_{x,t}(\lambda)$  (as  $ \Sdun \weak \Sdu$ in $ L^1(Q)^{3\times 3}$) so by  \eqref{g1} and \eqref{2g}   the RHS \eqref{rhslambda} is  non-positive as well as RHS of \eqref{3g}. That implies that
\begin{align}
\label{summ1}
(\bS(x,\vrn,\tn ,\Dun) - \bS(x,\vrn,\tn ,\Du) ):(\Dun - \Du) \to_b 0 .
\end{align}
What is more in the similar way as  \eqref{4g}  we obtain that 
\begin{align}
\label{summ2}
\bS(x,\vrn,\tn ,\Du):(\Dun - \Du) \to_b 0,  
\end{align}
and  one can infer also 
\begin{align}
\label{summ3}
\Sdun :\Du \to_b \Sdu :\Du. 
\end{align}
Summarising \eqref{summ1}--\eqref{summ3}  we provide  that  $ a_n \to_b a $  so assumption {\bf{A3}} of Lemma~\ref{bitlemma} holds. From the statement of Lemma \ref{bitlemma} we conclude that 
\begin{align*}
\Sdun :\Dun \weak \Sdu :\Du \weakly L^1 (Q) .
\end{align*}
This finishes the proof of Theorem \ref{main}.

\medskip
{\bf Acknowledgements: } 
The work of B.M. was supported by Grant of National Science Center Sonata, No 2013/09/D/ST1/03692 and his one-week stay in IMPAS was supported by WCMCS. A.W-K. was partially supported by Grant of National Science Center OPUS 4 2012/07/B/ST1/03306.



\end{document}